\documentclass[12pt,a4paper]{article}

\usepackage[utf8]{inputenc}             
\usepackage[english]{babel}

\usepackage{amsmath, amssymb, amsthm}
\usepackage{xcolor}
\usepackage{geometry}
\usepackage{graphicx,booktabs}
\usepackage{subcaption}  

\usepackage{authblk}  
\usepackage{hyperref}      
\usepackage{orcidlink} 

\usepackage{seqsplit}

\theoremstyle{plain}
\newtheorem{theorem}{Theorem}[section]
\newtheorem{lemma}[theorem]{Lemma}
\newtheorem{proposition}[theorem]{Proposition}
\newtheorem{corollary}[theorem]{Corollary}

\theoremstyle{definition}
\newtheorem{definition}[theorem]{Definition}

\newtheorem{example}[theorem]{Example}

\theoremstyle{remark}
\newtheorem{remark}[theorem]{Remark}

\title{Quantifying Dependence Effects on $k$-out-of-$n$ Reliability Relative to the Independent Benchmark}

\author[1,2]{
	Nuria Torrado\orcidlink{0000-0002-1156-3613} 
\thanks{nuria.torrado@gmail.com, nuria.torrado@uam.es, nuria.torrado@icmat.es}
}

\affil[1]{Department of Mathematics, Universidad Aut\'onoma de Madrid, 28049 Madrid, Spain}
\affil[2]{Institute of Mathematical Sciences (ICMAT), Campus de Cantoblanco, 28049 Madrid, Spain}

\date{\today}

\begin{document}
	\maketitle
	
%

\begin{abstract}
	Understanding how dependence among component lifetimes affects
	the reliability of redundant systems is a central problem in
	reliability engineering. In this paper, we develop a framework
	for quantifying the impact of multivariate dependence on
	$k$-out-of-$n$ system reliability relative to the independent
	benchmark. Using an orthogonal Legendre expansion of the copula
	density, we derive an exact decomposition that separates the
	dependence structure from the system architecture. This
	representation identifies how individual dependence contributions
	are transmitted to system reliability through explicit structural
	weights. 
	An application to multivariate
	Farlie--Gumbel--Morgenstern copulas demonstrates that higher-order
	dependence can produce distinct reliability effects even when
	all bivariate copulas remain unchanged. Our findings provide a
	systematic approach to quantifying and interpreting the effects
	of multivariate dependence on the reliability of redundant systems.
\end{abstract}

\noindent\textbf{Keywords:}
System reliability;
$k$-out-of-$n$ systems;
dependent components;
copula-based reliability analysis;
multivariate dependence;
higher-order dependence;
orthogonal Legendre expansions.


\section{Introduction}
\label{sec:introduction}

The $k$-out-of-$n$ structure is one of the basic redundancy models in
reliability theory. The system functions whenever at least $k$ of its
$n$ components are functioning, so that series and parallel systems are
recovered as the particular cases $k=n$ and $k=1$, respectively; see,
for example, \cite{rausand2004system}. 
%
%
 Different
aspects of the reliability of $k$-out-of-$n$ systems have been extensively
studied in the literature. In particular, asymptotic arguments can be used
to compare such systems under weak assumptions on the component lifetime
distributions \cite{torradoVeerman2012}. 
Reliability assessment has also been extended to weighted
$k$-out-of-$n$ systems with multiple types of components
\cite{hamdan2024weighted}.

When the component lifetimes are
independent, system reliability is determined by their marginal
distributions and the system structure. However, components operating
within the same system are often subject to common environmental
conditions, shared loads, interactions, or common failure mechanisms.
Their lifetimes may therefore be statistically dependent, and the
reliability analysis must take their joint distribution into account.
Copulas provide a natural framework for this purpose, since they separate
the marginal distributions from the dependence structure through Sklar's
theorem \cite{nelsen2006copulas}. This approach has been used in different
reliability settings involving dependent components. Copula-based
representations have been obtained for residual lifetimes of coherent
systems \cite{navarroDurante2017}, and relationships between properties of
the copula and reliability properties of coherent and $k$-out-of-$n$
systems have been investigated in \cite{navarroDuranteFernandez2021}.
Stochastic comparisons of conditional residual lifetimes
under different dependence structures, including FGM and
Archimedean copulas, have also been studied
\cite{zhang2023conditional}.
Distribution-free comparisons of coherent systems with heterogeneous and
dependent components have also been derived for general copula dependence
structures \cite{torrado2021coherent}. These contributions show that
dependence cannot in general be separated from the system structure when
the reliability of a coherent system is evaluated.

The analysis of dependent systems remains an active research topic.
Recent contributions include copula-based reliability models for
consecutive $k$-out-of-$n$ systems with different types of dependent
components \cite{ozkut2024consecutive}, series-parallel systems with
exchangeable dependent components \cite{torradoOzkut2024}, and
dependent-system reliability formulations combining order statistics,
survival signatures, and copula models \cite{leiKuo2025}. Copula methods
have also recently been applied to redundancy allocation problems for
series and parallel systems \cite{kumarRaoNaqvi2025}. These approaches
address different aspects of system reliability once a dependence model
has been specified and confirm the relevance of dependence modelling in
redundant systems.

A related problem arises in reliability sensitivity analysis when the
inputs are statistically dependent. In this setting, the interpretation
of sensitivity measures requires distinguishing the effects associated
with the probabilistic dependence model from those generated by the
system model itself \cite{ehrePapaioannouStraub2024}. A similar
distinction is relevant for $k$-out-of-$n$ systems. A specified copula,
or a summary measure such as a rank correlation coefficient, describes
dependence from the probabilistic side, whereas the effect of that
dependence on reliability also depends on the system architecture and on
the operating threshold. What is less explicit in the existing
copula-based reliability literature is a representation that separates
these two contributions at the level of individual components of the
multivariate dependence structure.

A natural question is therefore to quantify the change in system
reliability produced by dependence relative to the same system under
independent component lifetimes. For fixed marginal distributions, this
amounts to studying the difference between the reliability under the
actual dependence structure and the corresponding independence
benchmark. Beyond its overall magnitude, it is also relevant to identify
which features of the multivariate dependence structure are responsible
for this difference and how their effects vary with the system
architecture and the operating threshold.

Orthogonal representations of copula densities provide a suitable
framework for analysing multivariate dependence in a structured way.
Orthogonal series methods for copula densities have previously been
studied mainly from an estimation perspective; see
\cite{gui2009adaptive,ngounou2025orthproj}. Under square integrability,
the copula density can be represented in a tensor-product orthonormal
basis, so that its dependence structure is encoded through the
corresponding expansion coefficients. This representation is especially
convenient in the present reliability setting because it allows the
effect of dependence on the system reliability to be analysed term by
term. The present work establishes this connection for
$k$-out-of-$n$ systems and derives the system-side weights through which
the individual dependence contributions enter the reliability function.

In the present paper, we consider $k$-out-of-$n$ systems with dependent
components having identical continuous marginal lifetime distributions
and an absolutely continuous copula with square-integrable density.
The central question considered here is the reliability
difference between the dependent system and its independent counterpart
with the same marginal lifetime distribution. We refer to this difference
as the dependence impact. The orthogonal representation of the copula
density decomposes this impact into individual dependence contributions,
each multiplied by a structural weight determined by the system
architecture and the operating threshold.

The rest of the paper is organised as follows.
Section~\ref{sec:preliminaries} introduces the standardised
$k$-out-of-$n$ reliability function and recalls the orthogonal expansion
of the copula density in the shifted Legendre basis.
Section~\ref{sec:linear-form-weights} derives the linear decomposition of
the reliability function, introduces the dependence impact function, and
establishes the factorisation and the main structural properties of the
weights. Section~\ref{sec:fgm-applications} applies these results to the
multivariate FGM family and analyses the effect of an isolated
four-component dependence contribution while preserving the bivariate
copulas. Finally, Section~\ref{sec:conclusions} summarises the main
results and discusses possible extensions.

\section{Preliminaries}
\label{sec:preliminaries}

In this section, we recall the notation and the basic concepts that will be used throughout the paper. We first describe the $k$-out-of-$n$ reliability function on the original lifetime scale and then transform it to the uniform scale by means of the probability integral transform. After that, we state the orthogonal expansion of the copula density in the shifted Legendre basis. This representation is standard in orthogonal projection methods for copula densities, and we only retain the notation and identities needed in the sequel.

\subsection{Standardised $k$-out-of-$n$ reliability}
\label{subsec:system-model}

Let $\mathbf{T}=(T_1,\ldots,T_n)$ denote the vector of component lifetimes of a system with $n$ components. It is assumed that the marginal distributions are continuous and identical. We write
\[
F(t)=\mathbb{P}(T_i\le t),
\qquad
\bar F(t)=1-F(t),
\qquad i=1,\ldots,n,
\]
and denote by $\mathbf{F}$ the joint cumulative distribution function of $\mathbf{T}$.
If $T_{1:n}\le\cdots\le T_{n:n}$ are the order statistics associated with $\{T_i\}_{i=1}^n$, then the lifetime of the $k$-out-of-$n$ system is
$T_{k\mid n}:=T_{n-k+1:n}$.
Therefore, for $t\ge0$, the reliability function of the  $k$-out-of-$n$ system is given by
\begin{equation*}\label{eq:reliability-def-t}
	\bar{F}_{k\mid n}(t)
	:=
	\mathbb{P}(T_{k\mid n}>t)
	=
	\mathbb{P}\!\left(\sum_{i=1}^{n}\mathbb{I}\{T_i>t\}\ge k\right),
\end{equation*}
where $\mathbb{I}\{\cdot\}$ denotes the indicator function; see, for
instance, \cite{david2003order,rausand2004system}.
Equivalently, if
$X_i(t):=\mathbb{I}\{T_i>t\}$, $\mathbf{X}(t):=(X_1(t),\ldots,X_n(t))$ and
\[
\varphi_{k\mid n}(\mathbf{x})
:=
\mathbb{I}\Big\{\sum_{i=1}^n x_i\ge k\Big\},
\qquad \mathbf{x}\in\{0,1\}^n,
\]
then
$\bar{F}_{k\mid n}(t)=\mathbb{P}(\varphi_{k\mid n}(\mathbf{X}(t))=1)$.

By Sklar's theorem \cite{nelsen2006copulas}, there exists a unique copula
$\mathcal{C}:[0,1]^n\to[0,1]$ such that
\begin{equation*}\label{eq:sec2-sklar-T}
	\mathbf{F}(\mathbf{t})
	=
	\mathcal{C}\big(F(t_1),\ldots,F(t_n)\big),
	\qquad \mathbf{t}=(t_1,\ldots,t_n)\in\mathbb{R}^n.
\end{equation*}
Let us denote
\[
U_i:=F(T_i),
\qquad
\mathbf{U}=(U_1,\ldots,U_n).
\]
Then $U_i$ is uniformly distributed on $(0,1)$ for $i=1,\ldots,n$, and the copula $\mathcal{C}$ is the joint cumulative distribution function of $\mathbf{U}$, that is,
\begin{equation*}\label{eq:sec2-copula-is-law-of-U}
	\mathbb{P}(\mathbf{U}\le \mathbf{u})
	=
	\mathcal{C}(\mathbf{u}),
	\qquad \mathbf{u}=(u_1,\ldots,u_n)\in[0,1]^n,
\end{equation*}
where inequalities are understood componentwise. Thus, the transformation from $\mathbf{T}$ to $\mathbf{U}$ preserves the copula and removes the marginal distributions from the dependence analysis.

For $u=F(t)\in(0,1)$, the events $\{T_i>t\}$ and $\{U_i>u\}$ coincide. Hence, we define the standardised reliability function as
\begin{equation*}\label{eq:reliability-def-u}
	R_{k\mid n}(u)
	:=
	\mathbb{P}\!\left(\sum_{i=1}^{n}\mathbb{I}\{U_i>u\}\ge k\right),	
\end{equation*}
and thus, $\bar{F}_{k\mid n}(t)=R_{k\mid n}(F(t))$.
For $\mathbf{v}=(v_1,\ldots,v_n)\in[0,1]^n$, let
\begin{equation*}\label{eq:gku}
	g_{k,u}(\mathbf{v})
	:=
	\mathbb{I}\!\left\{\sum_{i=1}^{n}\mathbb{I}\{v_i>u\}\ge k\right\}
	=
	\varphi_{k\mid n}\big(\mathbb{I}\{v_1>u\},\ldots,\mathbb{I}\{v_n>u\}\big).
\end{equation*}
Therefore, $g_{k,u}$ describes the $k$-out-of-$n$ structure on the uniform scale. If the copula $\mathcal{C}$ is absolutely continuous with density $c$, then
\begin{equation}\label{eq:reliability-copula-integral}
	R_{k\mid n}(u)
	=
	\int_{[0,1]^n} g_{k,u}(\mathbf{v})\,c(\mathbf{v})\,\mathrm{d}\mathbf{v}
	=
	\mathbb{E}\big[g_{k,u}(\mathbf{U})\big].
\end{equation}
This identity separates the effect of the system structure, represented by $g_{k,u}$, from the dependence structure, represented by the copula density $c$.

In the independent case, $c(\mathbf{v})\equiv1$, and \eqref{eq:reliability-copula-integral} reduces to
\begin{equation}\label{eq:independence-benchmark}
	R_{k\mid n}^{\mathrm{ind}}(u)
	=
	\sum_{j=k}^{n}\binom{n}{j}(1-u)^{j}u^{n-j}.
\end{equation}
Throughout the paper, we use the conventions
\begin{equation}\label{eq:independence-benchmark-endpoint-conventions}
	R_{k\mid n}^{\mathrm{ind}}(u)=1 \quad \text{if } k\le 0,
	\qquad\textup{and}\qquad
	R_{k\mid n}^{\mathrm{ind}}(u)=0 \quad \text{if } k>n.
\end{equation}

\subsection{Copula density orthogonal expansion}
\label{subsec:legendre-notation}

Let $\{L_m\}_{m\ge0}$ be the classical Legendre polynomials on $[-1,1]$. Since the variables $U_i$ take values in $[0,1]$, we use the shifted orthogonal polynomials
$\mathbb{L}_m(x):=L_m(2x-1)$ with orthonormal polynomials
\begin{equation}\label{eq:bj-def}
	\phi_m(x):=\sqrt{2m+1}\:\mathbb{L}_m(x),
	\qquad m\ge0.
\end{equation}
Then $\phi_0(x)=1$ and
\begin{equation*}\label{eq:bj-orthonormal}
	\int_0^1\phi_j(x)\phi_m(x)\,\mathrm{d}x=\delta_{j,m},
	\qquad j,m\ge0.
\end{equation*}
For $\boldsymbol{m}=(m_1,\ldots,m_n)\in\mathbb{N}_0^n$, we define
\begin{equation}\label{eq:balpha-def}
	\Phi_{\boldsymbol{m}}(\mathbf{v})
	:=
	\prod_{i=1}^{n}\phi_{m_i}(v_i)
	\qquad\textup{and}\qquad
	\deg(\boldsymbol{m}):=m_1+\cdots+m_n.
\end{equation}
The family $\{\Phi_{\boldsymbol{m}}\}_{\boldsymbol{m}\in\mathbb{N}_0^n}$ is the usual tensor-product orthonormal basis of $L^2([0,1]^n)$. In particular,
\begin{equation*}\label{eq:balpha-orthonormal}
	\int_{[0,1]^n}\Phi_{\boldsymbol{m}}(\mathbf{u})\Phi_{\boldsymbol{n}}(\mathbf{u})\,\mathrm{d}\mathbf{u}
	=
	\delta_{\boldsymbol{m},\boldsymbol{n}},
	\qquad
	\boldsymbol{m},\boldsymbol{n}\in\mathbb{N}_0^n.
\end{equation*}
Classical properties of Legendre polynomials can be found in \cite{szego1975orthogonal}.

Assume that $\mathcal{C}$ is absolutely continuous and that its density $c$ satisfies
\begin{equation*}\label{eq:L2-assumption}
	\|c\|_{L^2([0,1]^n)}^2
	:=
	\int_{[0,1]^n}c(\mathbf{v})^2\,\mathrm{d}\mathbf{v}
	<\infty.
\end{equation*}
Under this assumption, Proposition~1 in \cite{ngounou2025orthproj} yields the expansion
\begin{equation*}\label{eq:sec3-copula-expansion}
	c(\mathbf{v})
	=
	\sum_{\boldsymbol{m}\in\mathbb{N}_0^n}
	\rho_{\boldsymbol{m}}\Phi_{\boldsymbol{m}}(\mathbf{v}),
\end{equation*}
where
\begin{equation}\label{eq:sec3-theta-def}
	\rho_{\boldsymbol{m}}
	:=
	\left\langle c,\Phi_{\boldsymbol{m}}\right\rangle_{L^2([0,1]^n)}
	=
	\int_{[0,1]^n}c(\mathbf{v})\Phi_{\boldsymbol{m}}(\mathbf{v})\,\mathrm{d}\mathbf{v}
	=
	\mathbb{E}\!\left[\Phi_{\boldsymbol{m}}(\mathbf{U})\right].
\end{equation}
The numbers $\rho_{\boldsymbol{m}}$ are the copula coefficients associated with the shifted Legendre basis. 
Since $\Phi_{\mathbf{0}}(\mathbf{v})=1$ and $c$ integrates to one, $\rho_{\mathbf{0}}=1$. Hence,
\begin{equation}\label{eq:sec3-copula-expansion-split}
	c(\mathbf{v})
	=
	1+
	\sum_{\boldsymbol{m}\neq\mathbf{0}}
	\rho_{\boldsymbol{m}}\Phi_{\boldsymbol{m}}(\mathbf{v}),
\end{equation}
where $\boldsymbol{0}=(0,\ldots,0)$.
Finally, the marginal uniformity of $U_i$ implies, for every $m\ge1$,
\begin{equation*}\label{eq:sec3-univariate-mean-zero}
	\mathbb{E}[\phi_m(U_i)]
	=
	\int_0^1\phi_m(u)\,\mathrm{d}u
	=0
	\qquad\textup{and}\qquad
	\mathbb{E}[\phi_m(U_i)^2]
	=
	\int_0^1\phi_m(u)^2\,\mathrm{d}u
	=1.
\end{equation*}
Consequently, if $\boldsymbol{m}\neq\mathbf{0}$ has only one nonzero component, then $\rho_{\boldsymbol{m}}=0$. Therefore, the nonconstant terms that may contribute to the copula density involve at least two nonzero coordinates.
For the product copula, $c(\mathbf{v})\equiv1$, and thus $\rho_{\mathbf{0}}=1$ and $\rho_{\boldsymbol{m}}=0$ for all $\boldsymbol{m}\neq\mathbf{0}$ (see \cite{gui2009adaptive,ngounou2025orthproj}).


\section{Decomposition of the standardised reliability}\label{sec:linear-form-weights}

The aim of this section is to connect the standardised reliability functional introduced in Subsection~\ref{subsec:system-model} with the orthogonal representation of the copula density recalled in Subsection~\ref{subsec:legendre-notation}.
To do this, we first define the active index set
\begin{equation*}\label{eq:active-set}
	\mathcal{A}(\boldsymbol{m})
	:=
	\{\, i\in\{1,\ldots,n\}:\ m_i>0 \,\},
\end{equation*}
for a given degree multi-index $\boldsymbol{m}=(m_1,\ldots,m_n)\in\mathbb{N}_0^{n}$.

The tensor-product polynomial defined in \eqref{eq:balpha-def} depends
nontrivially only on the coordinates indexed by
$\mathcal{A}(\boldsymbol{m})$, since $m_i=0$ implies
$\phi_{m_i}(x)=\phi_0(x)=1$. We therefore define
\[
s(\boldsymbol{m})
:=
\mathrm{card}\big(\mathcal{A}(\boldsymbol{m})\big)
\]
as the number of active coordinates associated with the multi-index
$\boldsymbol{m}$. Consequently,
\eqref{eq:sec3-copula-expansion-split} can be rewritten as
\begin{equation}\label{eq:sec3-restricted-expansion}
	c(\mathbf{v})
	=
	1+
	\sum_{\substack{\boldsymbol{m}\in\mathbb{N}_0^n\\
			s(\boldsymbol{m})\ge2}}
	\rho_{\boldsymbol{m}}\Phi_{\boldsymbol{m}}(\mathbf{v}).
\end{equation}

From a reliability standpoint, $s(\boldsymbol{m})$ determines the number
of component coordinates jointly involved in the dependence descriptor
indexed by $\boldsymbol{m}$. Hence, small values of
$s(\boldsymbol{m})$ correspond to lower-order dependence contributions,
such as pairwise effects, whereas larger values describe joint
contributions involving larger groups of components.


This interpretation connects directly with the $k$-out-of-$n$ structure.
At a fixed operating level $u$, system functioning requires at least
$k$ standardised component lifetimes to exceed the threshold. Accordingly,
the effect of a descriptor indexed by $\boldsymbol{m}$ on system
reliability depends both on its copula coefficient
$\rho_{\boldsymbol{m}}$ and on how the coordinates in
$\mathcal{A}(\boldsymbol{m})$ enter the $k$-out-of-$n$ success event.
The following result makes this separation explicit.

\begin{theorem}\label{prop:sec4-linear-form}
	Assume that the copula $\mathcal{C}$ is absolutely continuous with density $c\in L^2([0,1]^n)$ and let $u=F(t)\in(0,1)$. Then the standardised reliability \eqref{eq:reliability-copula-integral} can be written as
	\begin{equation}\label{eq:sec4-reliability-linear-form}
		R_{k\mid n}(u)
		=
		R_{k\mid n}^{\mathrm{ind}}(u)
		+
		\sum_{s(\boldsymbol{m})\ge 2}
		\rho_{\boldsymbol{m}}\,w_{\boldsymbol{m}}^{k\mid n}(u),
		\qquad\textup{for } \boldsymbol{m}\in\mathbb{N}_0^n,
	\end{equation}
	where $R_{k\mid n}^{\mathrm{ind}}(u)$ is the independence benchmark \eqref{eq:independence-benchmark}, $\{\rho_{\boldsymbol{m}}\}$ are the dependence descriptors \eqref{eq:sec3-theta-def}, and
	\begin{equation}\label{eq:sec4-weight-def}
		w_{\boldsymbol{m}}^{k\mid n}(u)
		:=
		\int_{[0,1]^n} g_{k,u}(\mathbf{v})\,\Phi_{\boldsymbol{m}}(\mathbf{v})\,\mathrm{d}\mathbf{v},
		\qquad \boldsymbol{m}\in\mathbb{N}_0^n.
	\end{equation}
	In particular, $w_{\mathbf{0}}^{k\mid n}(u)=R_{k\mid n}^{\mathrm{ind}}(u)$.
\end{theorem}

\begin{proof}
	By inserting the orthogonal expansion \eqref{eq:sec3-restricted-expansion} into \eqref{eq:reliability-copula-integral} and using linearity of the integral, we obtain
	\[
	R_{k\mid n}(u)
	=
	\int_{[0,1]^n} g_{k,u}(\mathbf{v})\,\mathrm{d}\mathbf{v}
	+
	\sum_{\substack{\boldsymbol{m}\in\mathbb{N}_0^n\\ s(\boldsymbol{m})\ge 2}}
	\rho_{\boldsymbol{m}}
	\int_{[0,1]^n} g_{k,u}(\mathbf{v})\,\Phi_{\boldsymbol{m}}(\mathbf{v})\,\mathrm{d}\mathbf{v}.
	\]
	The first integral is equal to $R_{k\mid n}^{\mathrm{ind}}(u)$ by \eqref{eq:independence-benchmark}, whereas the second one is $w_{\boldsymbol{m}}^{k\mid n}(u)$ by \eqref{eq:sec4-weight-def}. Finally, since $\Phi_{\mathbf{0}}(\mathbf{v})=1$,
	\[
	w_{\mathbf{0}}^{k\mid n}(u)
	=
	\int_{[0,1]^n} g_{k,u}(\mathbf{v})\,\mathrm{d}\mathbf{v}
	=
	R_{k\mid n}^{\mathrm{ind}}(u),
	\]
	and the proof is complete.
\end{proof}

The decomposition \eqref{eq:sec4-reliability-linear-form} separates the reliability of the system under independence from the part generated by dependence. The term $R_{k\mid n}^{\mathrm{ind}}(u)$ is the reliability of the same $k$-out-of-$n$ structure when the standardised variables $(U_1,\ldots,U_n)$ are independent, or equivalently when $c(\mathbf{u})\equiv1$ on $[0,1]^n$. Therefore it provides the natural benchmark level. The remaining series is exactly the deviation from this benchmark. This motivates the following definition.

\begin{definition}\label{def:dependence-impact}
	For a given $k$-out-of-$n$ system and a fixed marginal lifetime model, the
	\emph{dependence impact function} is defined as
	\begin{equation}\label{eq:dependence-impact}
		\mathcal{D}_{k\mid n}(u)
		:= R_{k\mid n}(u)-R^{\mathrm{ind}}_{k\mid n}(u),
		\qquad u\in(0,1).
	\end{equation}
\end{definition}

From the above definition and the bound
$0\le R_{k\mid n}(u)\le 1$, it follows that
\[
-\,R^{\mathrm{ind}}_{k\mid n}(u)
\le
\mathcal{D}_{k\mid n}(u)
\le
1-R^{\mathrm{ind}}_{k\mid n}(u),
\qquad u\in(0,1).
\]
Thus, $\mathcal{D}_{k\mid n}(u)$ measures the signed displacement of the reliability curve with respect to the independence case.
Positive values mean that dependence increases the system reliability at threshold $u$, whereas negative values mean that dependence decreases it. In particular, $\mathcal{D}_{k\mid n}(u)\le0$ for all $u\in(0,1)$ if and only if the dependent $k$-out-of-$n$ system is stochastically smaller than the independent benchmark, in the sense of the usual stochastic order. Analogously, $\mathcal{D}_{k\mid n}(u)\ge0$ for all $u\in(0,1)$ if and only if it is stochastically larger.
Moreover, from \eqref{eq:sec4-reliability-linear-form},
\begin{equation}\label{eq:dependence-impact-series}
	\mathcal{D}_{k\mid n}(u)=\sum_{\substack{\boldsymbol{m}\in\mathbb{N}_0^n\\ s(\boldsymbol{m})\ge 2}}\rho_{\boldsymbol{m}}\,w_{\boldsymbol{m}}^{k\mid n}(u).
\end{equation}
This expression shows how the departure from independence is formed. Each term contains a copula coefficient $\rho_{\boldsymbol{m}}$, which describes the dependence component indexed by $\boldsymbol{m}$, and a structural weight $w_{\boldsymbol{m}}^{k\mid n}(u)$, which indicates how the $k$-out-of-$n$ architecture transmits that component to the reliability functional at the threshold $u$. Consequently, a dependence descriptor may be large in the copula expansion and still have a small effect on reliability if its structural weight is small. Conversely, a moderate descriptor can be important for reliability when the corresponding structural weight is large.

\begin{remark}
	\label{rem:dependence-impact-differences}
	Consider two dependence specifications for the same $k$-out-of-$n$
	architecture and the same marginal lifetime model, arising, for example,
	from different copula families or from different parameter values within
	the same copula family. Since both specifications share the same
	independence benchmark $R_{k\mid n}^{\mathrm{ind}}(u)$, the Definition \ref{def:dependence-impact}
	gives
	\begin{equation*}
		\label{eq:dependence-impact-difference-identity}
		\mathcal{D}_{k\mid n}^{(2)}(u)-\mathcal{D}_{k\mid n}^{(1)}(u)
		=
		R_{k\mid n}^{(2)}(u)-R_{k\mid n}^{(1)}(u).
	\end{equation*}
	Moreover, by Theorem~\ref{prop:sec4-linear-form} and
	\eqref{eq:dependence-impact-series},
	\begin{equation}
		\label{eq:dependence-impact-coefficient-difference}
		R_{k\mid n}^{(2)}(u)-R_{k\mid n}^{(1)}(u)
		=
		\mathcal{D}_{k\mid n}^{(2)}(u)-\mathcal{D}_{k\mid n}^{(1)}(u)
		=
		\sum_{\substack{\boldsymbol{m}\in\mathbb{N}_0^n\\
				s(\boldsymbol{m})\geq2}}
		\left(
		\rho_{\boldsymbol{m}}^{(2)}
		-
		\rho_{\boldsymbol{m}}^{(1)}
		\right)
		w_{\boldsymbol{m}}^{k\mid n}(u).
	\end{equation}
	Hence, for a fixed system architecture, any change in reliability induced
	by a change in the dependence structure is exactly the corresponding
	change in dependence impact relative to independence. Equation
	\eqref{eq:dependence-impact-coefficient-difference} also identifies which
	changes in the copula coefficients are transmitted to the system and how
	they are weighted at the operating threshold $u$.
\end{remark}

The permutation symmetry of the $k$-out-of-$n$ structure induces a useful invariance property for the structural weights.
More precisely, for any permutation $\sigma$ of $\{1,\ldots,n\}$ and any multi-index $\boldsymbol{m}=(m_1,\ldots,m_n)\in\mathbb{N}_0^n$, write
\[
\sigma(\boldsymbol{m})
:=
(m_{\sigma(1)},\ldots,m_{\sigma(n)}).
\]
Since $g_{k,u}$ depends only on the number of coordinates above the threshold, it is invariant under relabelling of the components. Hence
\begin{equation}
	\label{eq:structural-weight-permutation-invariance}
	w_{\sigma(\boldsymbol{m})}^{k\mid n}(u)
	=
	w_{\boldsymbol{m}}^{k\mid n}(u),
	\qquad
	u\in(0,1).
\end{equation}
Therefore, for a fixed $k$-out-of-$n$ architecture, the structural contribution of a tensor-product term depends only on the multi-index pattern, up to permutation of component labels. This invariance is a property of the system structure and does not depend on the copula.

\subsection{Properties of the dependence impact function}
\label{subsec:dependence-impact-properties}

The function 
$\mathcal{D}_{k\mid n}(u)$, given by \eqref{eq:dependence-impact}, admits a natural interpretation as a reliability object.
Besides measuring the deviation from independence at a fixed threshold, it also satisfies several balance and symmetry properties. The following result shows that its integral has a direct meaning on the original lifetime scale.

\begin{proposition}
	\label{prop:mean-life-impact}
	Assume that the system lifetimes $T_{k\mid n}$ 
	and $T_{k\mid n}^{\mathrm{ind}}$ have finite means. Then the dependence impact function defined in \eqref{eq:dependence-impact} satisfies
	\begin{equation}
		\int_{0}^{\infty}\mathcal{D}_{k\mid n}(F(t))\,dt
		=
		\mathbb{E}[T_{k\mid n}]
		-
		\mathbb{E}[T_{k\mid n}^{\mathrm{ind}}],
		\label{eq:mean-life-impact-original-scale}
	\end{equation}
	where $T_{k\mid n}^{\mathrm{ind}}$ denotes the lifetime of the $k$-out-of-$n$ system under the independence benchmark with the same marginal distribution $F$.
\end{proposition}

\begin{proof}
	By \eqref{eq:dependence-impact}, for $u=F(t)\in(0,1)$,
	\begin{equation*}
		\mathcal{D}_{k\mid n}(F(t))
		=
		\mathbb{P}(T_{k\mid n}>t)
		-
		\mathbb{P}(T_{k\mid n}^{\mathrm{ind}}>t).
		\label{eq:proof-mean-life-step1}
	\end{equation*}
	Integrating over $t\in[0,\infty)$ gives
	\[
	\int_{0}^{\infty}\mathcal{D}_{k\mid n}(F(t))\,dt
	=
	\int_{0}^{\infty}\mathbb{P}(T_{k\mid n}>t)\,dt
	-
	\int_{0}^{\infty}\mathbb{P}(T_{k\mid n}^{\mathrm{ind}}>t)\,dt.
	\]
	Since the two lifetimes are nonnegative and have finite means, the tail-integral formula yields
	\[
	\int_{0}^{\infty}\mathbb{P}(T_{k\mid n}>t)\,dt
	=
	\mathbb{E}[T_{k\mid n}],
	\qquad
	\int_{0}^{\infty}\mathbb{P}(T_{k\mid n}^{\mathrm{ind}}>t)\,dt
	=
	\mathbb{E}[T_{k\mid n}^{\mathrm{ind}}],
	\]
	which proves \eqref{eq:mean-life-impact-original-scale}.
\end{proof}

For the standardised version, observe that $U_{k\mid n}:=F(T_{k\mid n})=U_{n-k+1:n}$ and $U_{k\mid n}^{\mathrm{ind}}:=F(T_{k\mid n}^{\mathrm{ind}})$. Then
\begin{equation*}\label{eq:mean-life-impact-standardised}
	\int_{0}^{1}\mathcal{D}_{k\mid n}(u)\,du
	=
	\mathbb{E}[U_{k\mid n}]
	-
	\mathbb{E}[U_{k\mid n}^{\mathrm{ind}}].	
\end{equation*}

Proposition~\ref{prop:mean-life-impact} shows that $\mathcal{D}_{k\mid n}$ is not only a local measure of the influence of dependence. Its integral gives the exact change in the mean lifetime of the system with respect to the independence benchmark. Hence the curve $\mathcal{D}_{k\mid n}(u)$ gives both a pointwise and a global description of the effect of dependence on system performance.

The next result gives a balance identity over the whole family of $k$-out-of-$n$ systems.

\begin{proposition}
	\label{prop:sum-D-zero}
	For every $u\in(0,1)$ and any $n\geq 2$, the dependence impact functions satisfy
	\begin{equation}
		\sum_{k=1}^{n}\mathcal{D}_{k\mid n}(u)=0.
		\label{eq:sum-D-zero}
	\end{equation}
\end{proposition}

\begin{proof}
	For each $k=1,\ldots,n$,
	\[
	R_{k\mid n}(u)
	=
	\mathbb{P}(U_{n-k+1:n}>u).
	\]
	Therefore,
	\[
	\sum_{k=1}^{n}R_{k\mid n}(u)
	=
	\sum_{j=1}^{n}\mathbb{P}(U_{j:n}>u)
	=
	\mathbb{E}\left[
	\sum_{j=1}^{n}\mathbb{I}_{\{U_{j:n}>u\}}
	\right].
	\]
	Ordering the observations does not change how many of them exceed a fixed threshold; see, e.g., \cite{david2003order}. Hence
	\[
	\sum_{j=1}^{n}\mathbb{I}_{\{U_{j:n}>u\}}
	=
	\sum_{i=1}^{n}\mathbb{I}_{\{U_i>u\}}.
	\]
	Taking expectations and using the uniform margins, we obtain
	\[
	\sum_{k=1}^{n}R_{k\mid n}(u)
	=
	\sum_{i=1}^{n}\mathbb{P}(U_i>u)
	=
	n(1-u).
	\]
	The same argument applied to the independent benchmark gives
	\[
	\sum_{k=1}^{n}R_{k\mid n}^{\mathrm{ind}}(u)
	=
	n(1-u).
	\]
	Consequently,
	\[
	\sum_{k=1}^{n}\mathcal{D}_{k\mid n}(u)
	=
	\sum_{k=1}^{n}R_{k\mid n}(u)
	-
	\sum_{k=1}^{n}R_{k\mid n}^{\mathrm{ind}}(u)
	=
	0,
	\]
	and \eqref{eq:sum-D-zero} follows.
\end{proof}

Proposition~\ref{prop:sum-D-zero} means that, at a fixed threshold, dependence redistributes reliability impact across the different $k$-out-of-$n$ architectures. Thus, if dependence produces a positive impact for some values of $k$, this must be compensated by a negative impact for other values. In this sense, the effect of dependence is inherently structure-dependent, as the same departure from independence cannot improve or degrade all $k$-out-of-$n$ systems simultaneously.

\begin{example}\label{ex:sum-D-zero-n3}
	Let us consider the case \(n=3\). For \(i<j\), denote by
	\[
	\mathcal{C}_{ij}(u,u)
	:=
	\mathbb{P}(U_i\le u,\,U_j\le u),
	\qquad u\in(0,1),
	\]
	the bivariate marginal copula evaluated on the diagonal, and write
	\[
	\mathcal{C}(u,u,u)
	=
	\mathbb{P}(U_1\le u,\,U_2\le u,\,U_3\le u).
	\]
	Then, the three standardised reliability functions are
	\begin{align*}
		R_{1\mid 3}(u)
		&=
		1-\mathcal{C}(u,u,u),\\
		R_{2\mid 3}(u)
		&=
		1-\mathcal{C}_{12}(u,u)-\mathcal{C}_{13}(u,u)-\mathcal{C}_{23}(u,u)
		+2\mathcal{C}(u,u,u),\\
		R_{3\mid 3}(u)
		&=
		1-3u+\mathcal{C}_{12}(u,u)+\mathcal{C}_{13}(u,u)
		+\mathcal{C}_{23}(u,u)-\mathcal{C}(u,u,u).
	\end{align*}
	Under independence, these functions reduce to
	\[
	R_{1\mid 3}^{\mathrm{ind}}(u)=1-u^3,
	\qquad
	R_{2\mid 3}^{\mathrm{ind}}(u)=1-3u^2+2u^3,
	\qquad
	R_{3\mid 3}^{\mathrm{ind}}(u)=(1-u)^3.
	\]
	Hence,
	\begin{align*}
		\mathcal{D}_{1\mid 3}(u)
		&=
		u^3-\mathcal{C}(u,u,u),\\
		\mathcal{D}_{2\mid 3}(u)
		&=
		3u^2-2u^3
		-\mathcal{C}_{12}(u,u)-\mathcal{C}_{13}(u,u)-\mathcal{C}_{23}(u,u)
		+2\mathcal{C}(u,u,u),\\
		\mathcal{D}_{3\mid 3}(u)
		&=
		\mathcal{C}_{12}(u,u)+\mathcal{C}_{13}(u,u)+\mathcal{C}_{23}(u,u)
		-\mathcal{C}(u,u,u)-3u^2+u^3.
	\end{align*}
	Adding the three expressions gives
	\[
	\mathcal{D}_{1\mid 3}(u)
	+
	\mathcal{D}_{2\mid 3}(u)
	+
	\mathcal{D}_{3\mid 3}(u)
	=
	0,
	\qquad u\in(0,1),
	\]
	in agreement with Proposition~\ref{prop:sum-D-zero}. This example also shows how the balance identity is achieved. The terms involving the trivariate copula and the bivariate marginal copulas do not disappear separately for each architecture, but they cancel when the three \(k\)-out-of-\(3\) configurations are considered together. Thus, dependence can move reliability impact from one architecture to another, but it cannot create a net impact over the whole family of \(k\)-out-of-\(3\) systems.
\end{example}

We now examine the case of exchangeable copulas. A further simplification is possible in this case, because terms that differ only by a permutation of component labels have the same copula coefficient and the same structural weight.

\begin{proposition}
	\label{prop:exchangeable-copula-grouping}
	Assume that the copula $\mathcal{C}$ is exchangeable and absolutely continuous with density $c\in L^2([0,1]^n)$. Then, for every permutation $\sigma$ of $\{1,\ldots,n\}$ and every multi-index $\boldsymbol{m}\in\mathbb{N}_0^n$,
	\begin{equation}
		\rho_{\sigma(\boldsymbol{m})}=\rho_{\boldsymbol{m}}.
		\label{eq:exchangeable-rho-invariance}
	\end{equation}
	Consequently, the dependence impact function can be regrouped over permutation classes as
	\begin{equation}
		\mathcal{D}_{k\mid n}(u)
		=
		\sum_{\substack{[\boldsymbol{m}]\\ s(\boldsymbol{m})\ge 2}}
		\mathrm{card}\bigl([\boldsymbol{m}]\bigr)\,
		\rho_{\boldsymbol{m}}\,
		w_{\boldsymbol{m}}^{k\mid n}(u),
		\qquad u\in(0,1),
		\label{eq:exchangeable-grouped-D}
	\end{equation}
	where
	$
	[\boldsymbol{m}]
	:=
	\{\sigma(\boldsymbol{m}):\sigma\in\mathfrak{S}_n\}
	$
	is the permutation class of $\boldsymbol{m}$, and the sum is taken over one representative from each class.
\end{proposition}

\begin{proof}
	Let
	$
	P_{\sigma}(\mathbf{v})
	:=
	(v_{\sigma(1)},\ldots,v_{\sigma(n)}),
	$
	for $\mathbf{v}\in[0,1]^n$. Since the copula is exchangeable, its density satisfies
	\[
	c\bigl(P_{\sigma}(\mathbf{v})\bigr)=c(\mathbf{v})
	\]
	for almost every $\mathbf{v}\in[0,1]^n$. By definition,
	\[
	\rho_{\sigma(\boldsymbol{m})}
	=
	\int_{[0,1]^n}
	c(\mathbf{v})\,
	\Phi_{\sigma(\boldsymbol{m})}(\mathbf{v})\,
	\mathrm{d}\mathbf{v}.
	\]
	Using the change of variables $\mathbf{v}=P_{\sigma}(\mathbf{x})$, and the fact that $P_{\sigma}$ is a measure preserving permutation of $[0,1]^n$, we get
	\[
	\rho_{\sigma(\boldsymbol{m})}
	=
	\int_{[0,1]^n}
	c\bigl(P_{\sigma}(\mathbf{x})\bigr)\,
	\Phi_{\sigma(\boldsymbol{m})}\bigl(P_{\sigma}(\mathbf{x})\bigr)\,
	\mathrm{d}\mathbf{x}.
	\]
	Moreover,
	\[
	\Phi_{\sigma(\boldsymbol{m})}\bigl(P_{\sigma}(\mathbf{x})\bigr)
	=
	\prod_{i=1}^n \phi_{m_{\sigma(i)}}(x_{\sigma(i)})
	=
	\prod_{i=1}^n \phi_{m_i}(x_i)
	=
	\Phi_{\boldsymbol{m}}(\mathbf{x}).
	\]
	Therefore,
	\[
	\rho_{\sigma(\boldsymbol{m})}
	=
	\int_{[0,1]^n}
	c(\mathbf{x})\,\Phi_{\boldsymbol{m}}(\mathbf{x})\,
	\mathrm{d}\mathbf{x}
	=
	\rho_{\boldsymbol{m}},
	\]
	which proves \eqref{eq:exchangeable-rho-invariance}. The grouped expression \eqref{eq:exchangeable-grouped-D} follows from \eqref{eq:dependence-impact-series}, from the  copula-coefficient invariance just proved, and from the permutation invariance of the structural weights \eqref{eq:structural-weight-permutation-invariance}. 
	Hence, all
	multi-indices in the same permutation class have the same coefficient
	and the same structural weight, so their total contribution is the
	common term multiplied by the cardinality of the class, $\mathrm{card}\bigl([\boldsymbol{m}]\bigr)$.
\end{proof}

We next consider reflection symmetry. Recall that a copula $\mathcal{C}$ is \textit{radially symmetric} if $\mathcal{C}=\widehat{\mathcal{C}}$, where $\widehat{\mathcal{C}}$ denotes the survival copula. For continuous copulas, this is equivalent to the distributional equality $(U_1,\ldots,U_n)\stackrel{d}{=}(1-U_1,\ldots,1-U_n)$; see \cite{nelsen2006copulas}. This property links the dependence impact functions of dual $k$-out-of-$n$ architectures through reflection of the threshold.

\begin{proposition}
	\label{prop:radial-duality-impact}
	Assume that the copula $\mathcal{C}$ is radially symmetric. Then, for every $k\in\{1,\ldots,n\}$ and every $u\in(0,1)$, the dependence impact function satisfies
	\begin{equation}\label{eq:radial-duality-D}
		\mathcal{D}_{k\mid n}(u)
		=
		-\mathcal{D}_{n-k+1\mid n}(1-u).		
	\end{equation}
\end{proposition}

\begin{proof}
	For the independent $k$-out-of-$n$ benchmark, the reflection identity for uniform order statistics yields
	\[
	R_{k\mid n}^{\mathrm{ind}}(1-u)
	=
	1-R_{n-k+1\mid n}^{\mathrm{ind}}(u),
	\qquad k=1,\ldots,n,
	\]
	where $U_{k:n}\stackrel{d}{=}1-U_{n-k+1:n}$ for order statistics of independent standard uniform random variables (see \cite{david2003order}). On the other hand, radial symmetry gives
	\[
	(1-U_1,\ldots,1-U_n)\stackrel{d}{=}(U_1,\ldots,U_n).
	\]
	Therefore, by the standard reflection property of order statistics,
	$
	U_{n-k+1:n}\stackrel{d}{=}1-U_{k:n}$,
	see \cite{david2003order,nelsen2006copulas}, and
	\[
	R_{k\mid n}(u)
	=
	\mathbb P(U_{n-k+1:n}>u)
	=
	1-\mathbb P(U_{k:n}>1-u)
	=
	1-R_{n-k+1\mid n}(1-u).
	\]
	Consequently, by Definition~\ref{def:dependence-impact}, we get
	\[
	\mathcal{D}_{k\mid n}(u)
	=
	R_{k\mid n}(u)-R_{k\mid n}^{\mathrm{ind}}(u)
	=
	\bigl[1-R_{n-k+1\mid n}(1-u)\bigr]
	-
	\bigl[1-R_{n-k+1\mid n}^{\mathrm{ind}}(1-u)\bigr]
	=
	-\mathcal{D}_{n-k+1\mid n}(1-u),
	\]
	and the result is proven.
\end{proof}

Proposition~\ref{prop:radial-duality-impact} states that, under radial symmetry, the dependence impact for a $k$-out-of-$n$ system at threshold $u$ is completely determined by the impact for the dual $(n-k+1)$-out-of-$n$ system at threshold $1-u$. The sign changes because survival of one dual structure is reflected into failure of the other. Thus, once the curve $\mathcal{D}_{k\mid n}(u)$ is known for a given value of $k$, the curve corresponding to its dual architecture is obtained directly from \eqref{eq:radial-duality-D} when the copula $\mathcal{C}$ is radially symmetric.

\begin{remark}
	\label{rem:radial-symmetry-rho}
	For the shifted Legendre basis used in this paper,
	the one dimensional basis
	polynomials satisfy $\phi_m(1-x)=(-1)^m\phi_m(x)$ (see \cite{szego1975orthogonal}). Hence
	\[
	\Phi_{\boldsymbol{m}}(\mathbf{1}-\mathbf{v})
	=
	(-1)^{\deg(\boldsymbol{m})}\Phi_{\boldsymbol{m}}(\mathbf{v}).
	\]
	If, in addition, the radially symmetric copula $\mathcal{C}$ admits a density $c$, then $c(\mathbf{u})=c(\mathbf{1}-\mathbf{u})$ for almost every $\mathbf{u}\in[0,1]^n$. From \eqref{eq:sec3-theta-def}, we obtain
	\begin{align*}
		\rho_{\boldsymbol{m}}
		&=
		\int_{[0,1]^n} c(\mathbf{u})\,\Phi_{\boldsymbol{m}}(\mathbf{u})\,d\mathbf{u}\\
		&=
		\int_{[0,1]^n} c(\mathbf{1}-\mathbf{u})\,\Phi_{\boldsymbol{m}}(\mathbf{u})\,d\mathbf{u}\\
		&=
		\int_{[0,1]^n} c(\mathbf{v})\,\Phi_{\boldsymbol{m}}(\mathbf{1}-\mathbf{v})\,d\mathbf{v}\\
		&=
		(-1)^{\deg(\boldsymbol{m})}
		\int_{[0,1]^n} c(\mathbf{v})\,\Phi_{\boldsymbol{m}}(\mathbf{v})\,d\mathbf{v}\\
		&=
		(-1)^{\deg(\boldsymbol{m})}\rho_{\boldsymbol{m}}.
	\end{align*}
	The third equality follows from the change of variables $\mathbf{v}=\mathbf{1}-\mathbf{u}$, which maps $[0,1]^n$ onto itself and has Jacobian determinant with absolute value equal to $1$. Therefore, $\rho_{\boldsymbol{m}}=0$ whenever $\deg(\boldsymbol{m})$ is odd. Thus, under radial symmetry, only even-degree dependence descriptors can appear in \eqref{eq:dependence-impact-series}.
\end{remark}

\subsection{Explicit representation and properties of the structural weights}\label{subsec:weights-computation}

For a fixed $k$-out-of-$n$ architecture, the weights $w_{\boldsymbol{m}}^{k\mid n}(u)$ do not depend on the copula. Hence, all copulas with square-integrable density share the same structural weights for a given system and threshold, only the coefficients $\rho_{\boldsymbol{m}}$ change. This subsection derives a finite expression for these weights and studies the resulting structural properties.

The definition \eqref{eq:sec4-weight-def} is an $n$-dimensional integral. However, the domain can be split according to the set of components whose standardised lifetimes exceed the operating threshold. For $A\subseteq\{1,\ldots,n\}$, let
\[
\mathcal{B}_A(u)
:=
\Big\{\mathbf{v}\in[0,1]^n:\ v_i>u \ \textup{for}\ i\in A,\ \ v_i\le u \ \textup{for}\ i\notin A\Big\}.
\]
On $\mathcal{B}_A(u)$, the indicator $g_{k,u}$ is equal to one if and only if $\operatorname{card}(A)\ge k$. The family $\{\mathcal{B}_A(u):A\subseteq\{1,\ldots,n\}\}$ partitions $[0,1]^n$ up to the null boundary where some coordinate is equal to $u$. Therefore
\[
w_{\boldsymbol{m}}^{k\mid n}(u)
=
\sum_{\substack{A\subseteq\{1,\ldots,n\}\\ \operatorname{card}(A)\ge k}}
\int_{\mathcal{B}_A(u)}\Phi_{\boldsymbol{m}}(\mathbf{v})\,\mathrm{d}\mathbf{v}.
\]
Using the product structure \eqref{eq:balpha-def}, each integral factorises into one-dimensional terms. For $m_i\ge0$, we define
\begin{equation*}\label{eq:sec4-Apm-def}
	A_{m_i}^{-}(u):=\int_{0}^{u}\phi_{m_i}(x)\,\mathrm{d}x,
	\qquad
	A_{m_i}^{+}(u):=\int_{u}^{1}\phi_{m_i}(x)\,\mathrm{d}x.
\end{equation*}
Then
\begin{equation}\label{eq:sec4-weight-block}
	w_{\boldsymbol{m}}^{k\mid n}(u)
	=
	\sum_{\substack{A\subseteq\{1,\ldots,n\}\\ \operatorname{card}(A)\ge k}}
	\prod_{i\in A} A_{m_i}^{+}(u)
	\prod_{i\notin A} A_{m_i}^{-}(u).
\end{equation}
This is already an exact finite representation. The system structure appears through $\operatorname{card}(A)\ge k$, while the basis appears through the one-dimensional primitives. Since $\phi_0(x)=1$,
\begin{equation}\label{eq:sec4-Apm}
	A_0^{-}(u)=u \quad\textup{and}\quad A_0^{+}(u)=1-u,
\end{equation}
and, by orthogonality with the constant polynomial, for every $m_i\ge1$,
\begin{equation}\label{eq:sec4-Apm-sign}
	A_{m_i}^{+}(u)=-A_{m_i}^{-}(u).
\end{equation}
The following theorem gives the compact factorisation used below.

\begin{theorem}
	\label{thm:structural-weight-factorisation}
	Let \(n\ge1\) and $k\in\{1,\ldots,n\}$. For every $\boldsymbol{m}\in\mathbb{N}_0^n$ and every $u\in(0,1)$, the structural weight defined in \eqref{eq:sec4-weight-def} satisfies
	\begin{equation}\label{eq:sec4-weight-reduction}
		w_{\boldsymbol{m}}^{k\mid n}(u)
		=
		\left(\prod_{i:m_i>0} A_{m_i}^{-}(u)\right)
		B_{s(\boldsymbol{m})}^{k\mid n}(u),
	\end{equation}
	where
	\begin{equation}\label{eq:Bkns-definition}
		B_{s}^{k\mid n}(u):=
		\sum_{j=0}^{s}
		\binom{s}{j}(-1)^j
		\,
		R_{k-j\mid n-s}^{\mathrm{ind}}(u),
		\qquad \textup{for }s=0,1,\ldots,n.
	\end{equation}
\end{theorem}

\begin{proof}
	We start from \eqref{eq:sec4-weight-block}. Fix a subset $A\subseteq\{1,\ldots,n\}$ and write
	\[
	j:=\operatorname{card}\bigl(A\cap\mathcal{A}(\boldsymbol{m})\bigr)
	\qquad\textup{and}\qquad
	\ell:=\operatorname{card}\bigl(A\cap\mathcal{A}(\boldsymbol{m})^{\mathrm c}\bigr).
	\]
	Thus $j$ is the number of active coordinates of $\boldsymbol{m}$ in $A$, and $\ell$ is the number of inactive coordinates in $A$. Since $\operatorname{card}(A)=j+\ell$, the condition $\operatorname{card}(A)\ge k$ is equivalent to $\ell\ge\max\{0,k-j\}$. Using \eqref{eq:sec4-Apm} for inactive coordinates and \eqref{eq:sec4-Apm-sign} for active coordinates, we have
	\begin{align*}
		\prod_{i\in A} A_{m_i}^{+}(u)
		\prod_{i\notin A} A_{m_i}^{-}(u)
		&=
		\left(
		\prod_{i\in A\cap \mathcal{A}(\boldsymbol{m})} A_{m_i}^{+}(u)
		\prod_{i\in \mathcal{A}(\boldsymbol{m})\setminus A} A_{m_i}^{-}(u)
		\right) \notag\\
		&\qquad{}
		\left(
		\prod_{i\in A\cap \mathcal{A}(\boldsymbol{m})^{\mathrm c}} A_0^{+}(u)
		\prod_{i\in \mathcal{A}(\boldsymbol{m})^{\mathrm c}\setminus A} A_0^{-}(u)
		\right) \notag\\
		&=
		\left(
		\prod_{i\in A\cap \mathcal{A}(\boldsymbol{m})} \bigl(-A_{m_i}^{-}(u)\bigr)
		\prod_{i\in \mathcal{A}(\boldsymbol{m})\setminus A} A_{m_i}^{-}(u)
		\right)
		(1-u)^{\ell}u^{\,n-s(\boldsymbol{m})-\ell} \notag\\
		&=
		(-1)^j
		\left(
		\prod_{i\in \mathcal{A}(\boldsymbol{m})} A_{m_i}^{-}(u)
		\right)
		(1-u)^{\ell}u^{\,n-s(\boldsymbol{m})-\ell}.
		\label{eq:app-product-fixedA}
	\end{align*}
	There are $\binom{s(\boldsymbol{m})}{j}$ ways to choose the active coordinates in $A$, and $\binom{n-s(\boldsymbol{m})}{\ell}$ ways to choose the inactive ones. Substitution in \eqref{eq:sec4-weight-block} gives
	\begin{align*}
		w_{\boldsymbol{m}}^{k\mid n}(u)
		&=
		\left(
		\prod_{i\in \mathcal{A}(\boldsymbol{m})} A_{m_i}^{-}(u)
		\right)
		\sum_{j=0}^{s(\boldsymbol{m})}
		\binom{s(\boldsymbol{m})}{j}(-1)^j
		\sum_{\ell=\max\{0,k-j\}}^{n-s(\boldsymbol{m})}
		\binom{n-s(\boldsymbol{m})}{\ell}
		(1-u)^{\ell}u^{\,n-s(\boldsymbol{m})-\ell}.
	\end{align*}
	By \eqref{eq:independence-benchmark}, the inner sum in the above equation is $R_{k-j\mid n-s(\boldsymbol{m})}^{\mathrm{ind}}(u)$. Hence
	\[
	w_{\boldsymbol{m}}^{k\mid n}(u)
	=
	\left(
	\prod_{i:m_i>0} A_{m_i}^{-}(u)
	\right)
	\sum_{j=0}^{s(\boldsymbol{m})}
	\binom{s(\boldsymbol{m})}{j}(-1)^j
	R_{k-j\mid n-s(\boldsymbol{m})}^{\mathrm{ind}}(u).
	\]
	Using \eqref{eq:Bkns-definition} with $s=s(\boldsymbol{m})$, this is exactly \eqref{eq:sec4-weight-reduction}.
\end{proof}

Theorem~\ref{thm:structural-weight-factorisation} separates two contributions. For a fixed multi-index $\boldsymbol{m}$ and a fixed threshold $u$, the product
\[
\prod_{i:m_i>0}A_{m_i}^{-}(u)
\]
is the same for all $k$-out-of-$n$ configurations. It only depends on the active coordinates and on their polynomial degrees. By contrast, $B_s^{k\mid n}(u)$ depends on the architecture $(k,n)$, the threshold $u$, and $s=s(\boldsymbol{m})$. Therefore, for the same dependence descriptor and the same threshold, the differences among the structural weights of the $k$-out-of-$n$ configurations are carried by the factors $B_s^{k\mid n}(u)$.

In reliability terms, $B_s^{k\mid n}(u)$ describes the contribution of the $n-s$ coordinates that remain after the active part has been isolated. If $j$ active coordinates are above the threshold, the remaining coordinates must provide at least $k-j$ additional exceedances for the system to work. This is why $R_{k-j\mid n-s}^{\mathrm{ind}}(u)$ appears in \eqref{eq:Bkns-definition}. It is important to note that the factor $B_s^{k\mid n}(u)$ is not itself a usual reliability function but it is a signed combination of independent reliability benchmarks. The signs come from \eqref{eq:sec4-Apm-sign}, which separates the exceedance and non-exceedance contributions of active coordinates.

\begin{example}\label{ex:structural-weights-11000}
	Let us consider a five-component system and the multi-index
	$
	\boldsymbol{m}=(1,1,0,0,0)$.
	Then
	$	\mathcal{A}(\boldsymbol{m})=\{1,2\}$, with
	$s(\boldsymbol{m})=2$.
	Since \(\phi_1(x)=\sqrt{3}(2x-1)\), we have
	\[
	A_1^{-}(u)
	=
	\int_0^u \phi_1(x)\,\mathrm{d}x
	=
	\sqrt{3}\int_0^u(2x-1)\,\mathrm{d}x
	=
	-\sqrt{3}\,u(1-u).
	\]
	Therefore,
	\[
	\prod_{i:m_i>0}A_{m_i}^{-}(u)
	=
	A_1^{-}(u)A_1^{-}(u)
	=
	3u^2(1-u)^2.
	\]	
	In this case, the factor in \eqref{eq:Bkns-definition} is
	\[
	B_2^{k\mid 5}(u)
	=
	R_{k\mid 3}^{\mathrm{ind}}(u)
	-
	2R_{k-1\mid 3}^{\mathrm{ind}}(u)
	+
	R_{k-2\mid 3}^{\mathrm{ind}}(u),
	\qquad k=1,\ldots,5.
	\]
	Using the endpoint conventions in
	\eqref{eq:independence-benchmark-endpoint-conventions}, together with
	\eqref{eq:independence-benchmark}, we have
	\[
	R_{1\mid 3}^{\mathrm{ind}}(u)=1-u^3,
	\qquad
	R_{2\mid 3}^{\mathrm{ind}}(u)=1-3u^2+2u^3,
	\qquad
	R_{3\mid 3}^{\mathrm{ind}}(u)=(1-u)^3,
	\]
	and therefore
	\begin{align*}
		B_2^{1\mid 5}(u)
		&=
		R_{1\mid 3}^{\mathrm{ind}}(u)
		-
		2R_{0\mid 3}^{\mathrm{ind}}(u)
		+
		R_{-1\mid 3}^{\mathrm{ind}}(u)
		=
		-u^3,\\
		B_2^{2\mid 5}(u)
		&=
		R_{2\mid 3}^{\mathrm{ind}}(u)
		-
		2R_{1\mid 3}^{\mathrm{ind}}(u)
		+
		R_{0\mid 3}^{\mathrm{ind}}(u)
		=
		u^2(4u-3),\\
		B_2^{3\mid 5}(u)
		&=
		R_{3\mid 3}^{\mathrm{ind}}(u)
		-
		2R_{2\mid 3}^{\mathrm{ind}}(u)
		+
		R_{1\mid 3}^{\mathrm{ind}}(u)
		=
		3u(1-u)(2u-1),\\
		B_2^{4\mid 5}(u)
		&=
		R_{4\mid 3}^{\mathrm{ind}}(u)
		-
		2R_{3\mid 3}^{\mathrm{ind}}(u)
		+
		R_{2\mid 3}^{\mathrm{ind}}(u)
		=
		(1-u)^2(4u-1),\\
		B_2^{5\mid 5}(u)
		&=
		R_{5\mid 3}^{\mathrm{ind}}(u)
		-
		2R_{4\mid 3}^{\mathrm{ind}}(u)
		+
		R_{3\mid 3}^{\mathrm{ind}}(u)
		=
		(1-u)^3.
	\end{align*}	
	Hence, by Theorem~\ref{thm:structural-weight-factorisation}, the five structural weights associated with
	\(\boldsymbol{m}=(1,1,0,0,0)\) are
	\begin{align*}
		w_{(1,1,0,0,0)}^{1\mid 5}(u)
		&=
		-3u^5(1-u)^2,\\
		w_{(1,1,0,0,0)}^{2\mid 5}(u)
		&=
		3u^4(1-u)^2(4u-3),\\
		w_{(1,1,0,0,0)}^{3\mid 5}(u)
		&=
		9u^3(1-u)^3(2u-1),\\
		w_{(1,1,0,0,0)}^{4\mid 5}(u)
		&=
		3u^2(1-u)^4(4u-1),\\
		w_{(1,1,0,0,0)}^{5\mid 5}(u)
		&=
		3u^2(1-u)^5.
	\end{align*}
	
	These expressions are independent of the copula. They depend only on the \(k\)-out-of-\(5\) architecture, the threshold \(u\), and the selected tensor-product term. The copula affects the reliability expansion only through the coefficient \(\rho_{(1,1,0,0,0)}\). Thus, for this pairwise dependence descriptor, the structural weights above indicate where and with which sign the system architecture can transmit the descriptor into the dependence impact function \(\mathcal{D}_{k\mid 5}(u)\).
	
	It is also worth noting that the structural weights obtained above differ across the five system configurations. This illustrates how the same dependence descriptor can have different structural effects depending on the system architecture.
\end{example}


The factorization \eqref{eq:sec4-weight-reduction} also provides the following boundary property.

\begin{proposition}\label{prop:wm-extremes}
	Let \(n\ge1\). Fix $k\in\{1,\ldots,n\}$ and let $\boldsymbol{m}\neq\boldsymbol{0}$. Then there exist constants $C_{\boldsymbol{m},k}^{(0)}>0$ and $C_{\boldsymbol{m},k}^{(1)}>0$ such that, for all $u\in(0,1)$,
	\begin{equation}\label{eq:wm-extremes-lower-upper}
		\left|w_{\boldsymbol{m}}^{k\mid n}(u)\right|
		\le
		C_{\boldsymbol{m},k}^{(0)}u^{s(\boldsymbol{m})},
		\qquad
		\left|w_{\boldsymbol{m}}^{k\mid n}(u)\right|
		\le
		C_{\boldsymbol{m},k}^{(1)}(1-u)^{s(\boldsymbol{m})}.
	\end{equation}
	Consequently,
	\[
	\lim_{u\downarrow0}w_{\boldsymbol{m}}^{k\mid n}(u)=0,
	\qquad
	\lim_{u\uparrow1}w_{\boldsymbol{m}}^{k\mid n}(u)=0.
	\]
\end{proposition}

\begin{proof}
	The factor $B_{s(\boldsymbol{m})}^{k\mid n}$ is a polynomial in $u$, and therefore it is bounded on $[0,1]$. Hence there exists $M_{\boldsymbol{m},k}>0$ such that
	\[
	\left|B_{s(\boldsymbol{m})}^{k\mid n}(u)\right|
	\le
	M_{\boldsymbol{m},k},
	\qquad 0\le u\le1.
	\]
	Let $i$ be an active index, so that $m_i\ge1$. By definition, $A^-_{m_i}(0)=0$, and by orthogonality with the constant basis function,
	\[
	A^-_{m_i}(1)
	=
	\int_0^1\phi_{m_i}(x)\,\mathrm{d}x
	=
	0.
	\]
	Since $A^-_{m_i}$ is a polynomial with zeros at $0$ and $1$, there exists a polynomial $Q_i$ such that
	\[
	A^-_{m_i}(u)=u(1-u)Q_i(u).
	\]
	Consequently, there exist constants $M_i^{(0)},M_i^{(1)}>0$ such that, for all $u\in[0,1]$,
	\[
	|A^-_{m_i}(u)|\le M_i^{(0)}u,
	\qquad
	|A^-_{m_i}(u)|\le M_i^{(1)}(1-u).
	\]
	Multiplying these bounds over the active indices and using \eqref{eq:sec4-weight-reduction}, we obtain
	\[
	\left|w_{\boldsymbol{m}}^{k\mid n}(u)\right|
	\le
	M_{\boldsymbol{m},k}
	\left(
	\prod_{i:m_i>0}M_i^{(0)}
	\right)
	u^{s(\boldsymbol{m})},
	\]
	and
	\[
	\left|w_{\boldsymbol{m}}^{k\mid n}(u)\right|
	\le
	M_{\boldsymbol{m},k}
	\left(
	\prod_{i:m_i>0}M_i^{(1)}
	\right)
	(1-u)^{s(\boldsymbol{m})}.
	\]
	Thus \eqref{eq:wm-extremes-lower-upper} follows by taking
	\[
	C_{\boldsymbol{m},k}^{(0)}
	=
	M_{\boldsymbol{m},k}
	\prod_{i:m_i>0}M_i^{(0)},
	\qquad
	C_{\boldsymbol{m},k}^{(1)}
	=
	M_{\boldsymbol{m},k}
	\prod_{i:m_i>0}M_i^{(1)}.
	\]
	The two limits are immediate from \eqref{eq:wm-extremes-lower-upper}.
\end{proof}

Proposition~\ref{prop:wm-extremes} has a direct operational interpretation. When $u$ is close to zero, the system is in an early-life regime on the marginal scale and most components exceed the threshold. When $u$ is close to one, the system is in a late-life regime and threshold exceedances are rare. In both regimes, the scope for dependency to modify reliability is limited, and the structural weights associated with nonconstant dependence descriptors vanish. The order of this attenuation is governed by $s(\boldsymbol{m})$, while the architecture-specific part is contained in $B_{s(\boldsymbol{m})}^{k\mid n}$. Thus dependence is expected to be most relevant in intermediate regimes, where the system is neither almost surely functioning nor almost surely failed.

\begin{remark}\label{rem:best-boundary-constants}
	The constants in Proposition~\ref{prop:wm-extremes} provide admissible
	constants for the upper bounds. A more precise characterisation of the
	boundary behaviour is obtained by considering the corresponding optimal
	admissible constants given by
	\begin{equation}\label{eq:optimal-boundary-constants}
		C_{\boldsymbol{m},k,\mathrm{opt}}^{(0)}
		:=
		\sup_{0<u<1}
		\frac{\left|w_{\boldsymbol{m}}^{k\mid n}(u)\right|}
		{u^{s(\boldsymbol{m})}},
		\qquad
		C_{\boldsymbol{m},k,\mathrm{opt}}^{(1)}
		:=
		\sup_{0<u<1}
		\frac{\left|w_{\boldsymbol{m}}^{k\mid n}(u)\right|}
		{(1-u)^{s(\boldsymbol{m})}}.
	\end{equation}
	Proposition~\ref{prop:wm-extremes} ensures that these quantities are
	finite. They are the smallest nonnegative constants that make the
	corresponding boundary bounds valid on $(0,1)$. Indeed, if a nonnegative
	constant $C$ satisfies
	\[
	\left|w_{\boldsymbol{m}}^{k\mid n}(u)\right|
	\le
	Cu^{s(\boldsymbol{m})},
	\qquad 0<u<1,
	\]
	then
	\[
	\frac{\left|w_{\boldsymbol{m}}^{k\mid n}(u)\right|}
	{u^{s(\boldsymbol{m})}}
	\le C,
	\qquad 0<u<1,
	\]
	and hence $C_{\boldsymbol{m},k,\mathrm{opt}}^{(0)}\le C$. Therefore,
	no smaller nonnegative constant can control the lower-threshold
	attenuation in the first bound of
	\eqref{eq:wm-extremes-lower-upper}. The same reasoning applies to
	$C_{\boldsymbol{m},k,\mathrm{opt}}^{(1)}$ for the upper-threshold
	attenuation. In reliability terms, these constants quantify the largest
	normalised structural contribution of the descriptor $\boldsymbol{m}$
	relative to the lower- and upper-threshold attenuation scales,
	respectively.
\end{remark}

We now study the structural duality obtained by replacing $u$ with $1-u$ and the $k$-out-of-$n$ architecture with its dual $(n-k+1)$-out-of-$n$ architecture. This property concerns the weights themselves and does not require any additional assumption on the copula.

\begin{lemma}\label{lem:Bkns-symmetry}
	Let $n\ge1$, $k\in\{1,\ldots,n\}$, and $s\in\{1,\ldots,n\}$. The factor $B_s^{k\mid n}(u)$ defined in \eqref{eq:Bkns-definition} satisfies
	\begin{equation}\label{eq:Bkns-symmetry}
		B_s^{k\mid n}(1-u)
		=
		(-1)^{s+1}B_s^{n-k+1\mid n}(u),
		\qquad u\in(0,1).
	\end{equation}
	In particular,
	\begin{equation}\label{eq:Bkns-extremes}
		B_s^{1\mid n}(u)=-u^{n-s}
		\qquad\textup{and}\qquad
		B_s^{n\mid n}(u)=(-1)^s(1-u)^{n-s}.
	\end{equation}
\end{lemma}

\begin{proof}
	For the independent $q$-out-of-$p$ benchmark, the reflection identity for uniform order statistics gives
	\[
	R_{q\mid p}^{\mathrm{ind}}(1-u)
	=
	1-R_{p-q+1\mid p}^{\mathrm{ind}}(u),
	\qquad q=1,\ldots,p,
	\]
	where $U_{q:p}\stackrel{d}{=}1-U_{p-q+1:p}$ (see \cite{david2003order}). By the endpoint convention in \eqref{eq:independence-benchmark-endpoint-conventions}, the same identity also covers the boundary cases that appear below. Applying it with $p=n-s$ and $q=k-j$, and using \eqref{eq:Bkns-definition}, we obtain
	\begin{align*}
		B_s^{k\mid n}(1-u)
		&=
		\sum_{j=0}^{s}
		\binom{s}{j}(-1)^j
		\left[
		1-
		R_{n-s-k+j+1\mid n-s}^{\mathrm{ind}}(u)
		\right].
	\end{align*}
	Since
	\[
	\sum_{j=0}^{s}\binom{s}{j}(-1)^j=0,
	\]
	the constant term cancels, and hence
	\[
	B_s^{k\mid n}(1-u)
	=
	-
	\sum_{j=0}^{s}
	\binom{s}{j}(-1)^j
	R_{n-s-k+j+1\mid n-s}^{\mathrm{ind}}(u).
	\]
	With the change of index $r=s-j$, this becomes
	\[
	B_s^{k\mid n}(1-u)
	=
	(-1)^{s+1}
	\sum_{r=0}^{s}
	\binom{s}{r}(-1)^r
	R_{n-k+1-r\mid n-s}^{\mathrm{ind}}(u).
	\]
	The last sum is $B_s^{n-k+1\mid n}(u)$, and this proves \eqref{eq:Bkns-symmetry}.
	
	For $k=1$, we have $R_{1\mid n-s}^{\mathrm{ind}}(u)=1-u^{n-s}$ and
	$R_{1-j\mid n-s}^{\mathrm{ind}}(u)=1$ for $j\ge 1$. Together with
	\eqref{eq:Bkns-definition}, this gives
	\begin{eqnarray*}		
		B_s^{1\mid n}(u)
		&=&
		R_{1\mid n-s}^{\mathrm{ind}}(u)
		+
		\sum_{j=1}^{s}
		\binom{s}{j}(-1)^j
		R_{1-j\mid n-s}^{\mathrm{ind}}(u)\\
		&=&
		1-u^{n-s}
		+
		\sum_{j=1}^{s}\binom{s}{j}(-1)^j
		=
		-u^{n-s}.
	\end{eqnarray*}	
	For $k=n$, all terms vanish except the one with $j=s$, and therefore
	\[
	B_s^{n\mid n}(u)
	=
	(-1)^sR_{n-s\mid n-s}^{\mathrm{ind}}(u)
	=
	(-1)^s(1-u)^{n-s}.
	\]
	This proves \eqref{eq:Bkns-extremes}.
\end{proof}

The preceding symmetry of $B_s^{k\mid n}(u)$, together with the parity of shifted Legendre polynomials, gives the following duality for the  structural weights.

\begin{proposition}\label{prop:dual_weights}
	Let \(n\ge1\) and $k\in\{1,\ldots,n\}$. For every multi-index $\boldsymbol{m}\in\mathbb{N}_0^n$ such that $\boldsymbol{m}\neq\boldsymbol{0}$ and every $u\in(0,1)$,
	\begin{equation}\label{eq:dual_weight_relation}
		w_{\boldsymbol{m}}^{k\mid n}(1-u)
		=
		(-1)^{1+\deg(\boldsymbol{m})}
		w_{\boldsymbol{m}}^{n-k+1\mid n}(u),
	\end{equation}
	where $\deg(\boldsymbol{m})=m_1+\cdots+m_n$. Consequently,
	\begin{equation}\label{eq:dual_D_relation}
		\mathcal{D}_{k\mid n}(1-u)
		=
		\sum_{\substack{\boldsymbol{m}\in\mathbb{N}_0^n\\ s(\boldsymbol{m})\ge2}}
		(-1)^{1+\deg(\boldsymbol{m})}\rho_{\boldsymbol{m}}
		w_{\boldsymbol{m}}^{n-k+1\mid n}(u).
	\end{equation}
\end{proposition}

\begin{proof}
	For $m_i\ge1$, the parity of shifted Legendre polynomials gives
	$
	\phi_{m_i}(1-x)=(-1)^{m_i}\phi_{m_i}(x)$,
	see, e.g., \cite{szego1975orthogonal}. Hence, using \eqref{eq:sec4-Apm-sign} and the change of variable $x=1-y$, we get
	\[
	A^-_{m_i}(1-u)
	=
	\int_0^{1-u}\phi_{m_i}(x)\,\mathrm{d}x
	=
	(-1)^{m_i}\int_u^1\phi_{m_i}(y)\,\mathrm{d}y
	=
	(-1)^{m_i+1}A^-_{m_i}(u).
	\]
	Multiplying over the active indices gives
	\[
	\prod_{i:m_i>0}A^-_{m_i}(1-u)
	=
	(-1)^{\deg(\boldsymbol{m})+s(\boldsymbol{m})}
	\prod_{i:m_i>0}A^-_{m_i}(u).
	\]
	Combining this identity with Lemma~\ref{lem:Bkns-symmetry} in \eqref{eq:sec4-weight-reduction} yields \eqref{eq:dual_weight_relation}, since $(-1)^{2s(\boldsymbol{m})}=1$. Substitution in \eqref{eq:dependence-impact-series} gives \eqref{eq:dual_D_relation}.
\end{proof}

\begin{example}\label{ex:dual-weights-11000}
	Let us continue with the multi-index considered in
	Example~\ref{ex:structural-weights-11000}, namely
	$\boldsymbol{m}=(1,1,0,0,0)$.
	In this case,
	$	\deg(\boldsymbol{m})=2$, and
	$s(\boldsymbol{m})=2$.
	Therefore, Proposition~\ref{prop:dual_weights} gives
	\[
	w_{\boldsymbol{m}}^{k\mid 5}(1-u)
	=
	-
	w_{\boldsymbol{m}}^{6-k\mid 5}(u),
	\qquad k=1,\ldots,5.
	\]
	This identity can be checked directly from the explicit expressions obtained in
	Example~\ref{ex:structural-weights-11000}.	
	Thus, for \(k=1\),
	\[
	w_{(1,1,0,0,0)}^{1\mid 5}(1-u)
	=
	-3(1-u)^5u^2
	=-w_{(1,1,0,0,0)}^{5\mid 5}(u).
	\]	
	Analogously, for \(k=2\),
	\[
	w_{(1,1,0,0,0)}^{2\mid 5}(1-u)
	=
	3(1-u)^4u^2\bigl(4(1-u)-3\bigr)
	=
	3u^2(1-u)^4(1-4u)
	=	-w_{(1,1,0,0,0)}^{4\mid 5}(u).
	\]	
	The remaining case \(k=3\) corresponds to the self-dual architecture. Since \(6-k=3\), the same identity becomes
	\begin{align*}
		w_{(1,1,0,0,0)}^{3\mid 5}(1-u)
		&=
		9(1-u)^3u^3\bigl(2(1-u)-1\bigr)\\
		&=
		9u^3(1-u)^3(1-2u)\\
		&=
		-
		9u^3(1-u)^3(2u-1)\\
		&=
		-
		w_{(1,1,0,0,0)}^{3\mid 5}(u).
	\end{align*}	
	Therefore, the duality relation in \eqref{eq:dual_weight_relation} is verified explicitly for the pairwise descriptor \((1,1,0,0,0)\). In reliability terms, the structural contribution assigned by the \(k\)-out-of-\(5\) system at the reflected threshold \(1-u\) is the negative of the contribution assigned by the dual \((6-k)\)-out-of-\(5\) system at threshold \(u\). Since this descriptor has even total degree, the sign in \eqref{eq:dual_weight_relation} is negative. Hence, the two dual architectures have structural weights with the same magnitude and opposite sign at reflected threshold levels.
\end{example}

As a consequence of Proposition~\ref{prop:dual_weights}, dual architectures have structural weights with the same absolute value at reflected thresholds.

\begin{corollary}\label{cor:dual_weight_magnitude_symmetry}
	Let \(n\ge1\) and $k\in\{1,\ldots,n\}$. For every multi-index $\boldsymbol{m}\in\mathbb{N}_0^n$ such that $\boldsymbol{m}\neq\boldsymbol{0}$, and every $u\in(0,1)$,
	\begin{equation}\label{eq:dual_weight_magnitude_symmetry}
		\left|w_{\boldsymbol{m}}^{k\mid n}(1-u)\right|
		=
		\left|w_{\boldsymbol{m}}^{n-k+1\mid n}(u)\right|.
	\end{equation}
	Moreover,
	\begin{equation}\label{eq:optimal-boundary-constants-duality}
		C_{\boldsymbol{m},k,\mathrm{opt}}^{(0)}
		=
		C_{\boldsymbol{m},\,n-k+1,\mathrm{opt}}^{(1)},
		\qquad
		C_{\boldsymbol{m},k,\mathrm{opt}}^{(1)}
		=
		C_{\boldsymbol{m},\,n-k+1,\mathrm{opt}}^{(0)}.
	\end{equation}
\end{corollary}

\begin{proof}
	The identity \eqref{eq:dual_weight_magnitude_symmetry} follows by taking absolute values in \eqref{eq:dual_weight_relation}. Using this identity in \eqref{eq:optimal-boundary-constants} and applying the change of variable $v=1-u$, we obtain
	\[
	C_{\boldsymbol{m},k,\mathrm{opt}}^{(0)}
	=
	\sup_{0<v<1}
	\frac{\left|w_{\boldsymbol{m}}^{n-k+1\mid n}(v)\right|}
	{(1-v)^{s(\boldsymbol{m})}}
	=
	C_{\boldsymbol{m},\,n-k+1,\mathrm{opt}}^{(1)}.
	\]
	The second identity in \eqref{eq:optimal-boundary-constants-duality} follows in the same way.
\end{proof}

Corollary~\ref{cor:dual_weight_magnitude_symmetry} shows that a $k$-out-of-$n$ system at threshold $1-u$ and its dual $(n-k+1)$-out-of-$n$ system at threshold $u$ assign the same magnitude to the structural contribution associated with the same descriptor. The identities in \eqref{eq:optimal-boundary-constants-duality} give the corresponding statement for the boundary bounds, in the sense that the optimal lower-threshold constant of one configuration is the optimal upper-threshold constant of its dual, and conversely.

In the self-dual case, namely when $n$ is odd and $k=(n+1)/2$, the dual system coincides with the original one. Therefore,
\[
\bigl|w_{\boldsymbol{m}}^{k\mid n}(1-u)\bigr|
=
\bigl|w_{\boldsymbol{m}}^{k\mid n}(u)\bigr|,
\qquad u\in(0,1).
\]
Moreover, the two optimal boundary constants coincide,
\[
C_{\boldsymbol{m},k,\mathrm{opt}}^{(0)}
=
C_{\boldsymbol{m},k,\mathrm{opt}}^{(1)}.
\]

The previous result compares reflected thresholds. A stronger pointwise relation is obtained for fully active descriptors, that is, for descriptors involving all components.

\begin{proposition}\label{prop:fully-active-dual-weights}
	Let $\boldsymbol{m}\in\mathbb{N}_0^n$ be such that $s(\boldsymbol{m})=n$. Let $n\geq 1$.Then, for every $k\in\{1,\ldots,n\}$ and every $u\in(0,1)$,
	\begin{equation}\label{eq:fully-active-dual-weight-relation}
		w_{\boldsymbol{m}}^{n-k+1\mid n}(u)
		=
		(-1)^{n+1}
		w_{\boldsymbol{m}}^{k\mid n}(u).
	\end{equation}
	Consequently,
	\begin{equation}\label{eq:fully-active-dual-weight-magnitude-relation}
		\left|w_{\boldsymbol{m}}^{n-k+1\mid n}(u)\right|
		=
		\left|w_{\boldsymbol{m}}^{k\mid n}(u)\right|.
	\end{equation}
\end{proposition}

\begin{proof}
	Since $s(\boldsymbol{m})=n$, all coordinates are active. Thus, by \eqref{eq:sec4-weight-reduction},
	\[
	w_{\boldsymbol{m}}^{k\mid n}(u)
	=
	\left(
	\prod_{i=1}^{n}A_{m_i}^{-}(u)
	\right)
	B_n^{k\mid n}(u).
	\]
	It remains to compare $B_n^{k\mid n}(u)$ and $B_n^{n-k+1\mid n}(u)$. From \eqref{eq:Bkns-definition},
	\[
	B_n^{k\mid n}(u)
	=
	\sum_{j=0}^{n}
	\binom{n}{j}(-1)^j
	R_{k-j\mid 0}^{\mathrm{ind}}(u).
	\]
	By the endpoint convention \eqref{eq:independence-benchmark-endpoint-conventions}, applied with $n=0$,
	\[
	R_{k-j\mid 0}^{\mathrm{ind}}(u)
	=
	\begin{cases}
		1, & j\ge k,\\
		0, & j<k.
	\end{cases}
	\]
	Therefore,
	\[
	B_n^{k\mid n}(u)
	=
	\sum_{j=k}^{n}\binom{n}{j}(-1)^j
	=
	(-1)^k\binom{n-1}{k-1},
	\]
	where the last equality is the elementary alternating binomial identity. Similarly,
	\[
	B_n^{n-k+1\mid n}(u)
	=
	\sum_{j=n-k+1}^{n}\binom{n}{j}(-1)^j
	=
	(-1)^{n-k+1}\binom{n-1}{n-k}.
	\]
	Since $\binom{n-1}{n-k}=\binom{n-1}{k-1}$, it follows that
	\[
	B_n^{n-k+1\mid n}(u)
	=
	(-1)^{n+1}B_n^{k\mid n}(u).
	\]
	The product $\prod_{i=1}^{n}A_{m_i}^{-}(u)$ is the same for the two dual architectures. Substituting the preceding identity into \eqref{eq:sec4-weight-reduction} gives \eqref{eq:fully-active-dual-weight-relation}. Taking absolute values gives \eqref{eq:fully-active-dual-weight-magnitude-relation}.
\end{proof}

Proposition~\ref{prop:fully-active-dual-weights} has a direct reliability interpretation. When a dependence descriptor involves all $n$ components, the dual $k$-out-of-$n$ and $(n-k+1)$-out-of-$n$ architectures assign the same structural magnitude to that descriptor at the same operating threshold. Thus \eqref{eq:fully-active-dual-weight-magnitude-relation} is not a reflected-threshold identity, but a pointwise comparison between dual architectures. If $n$ is odd, the structural weights coincide for all $u\in(0,1)$, whereas if $n$ is even, they only differ by a sign.

Finally, we close this subsection by displaying three representative structural weight profiles for a five component system.

\begin{figure}[ht]
	\centering
	\captionsetup[subfigure]{labelformat=simple,labelsep=space}
	\renewcommand\thesubfigure{\alph{subfigure})}
	
	\begin{subfigure}[t]{0.33\textwidth}
		\centering
		\includegraphics[width=\linewidth]{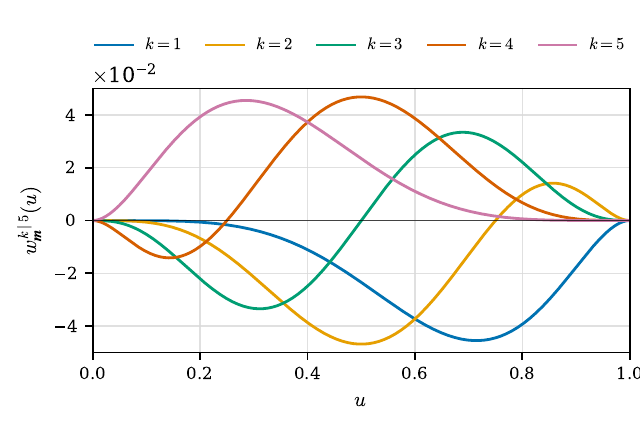}
		\caption{$\boldsymbol{m}_1=(1,1,0,0,0)$}
		\label{fig:structural-weights-11000-n5}
	\end{subfigure}\hfill
	\begin{subfigure}[t]{0.33\textwidth}
		\centering
		\includegraphics[width=\linewidth]{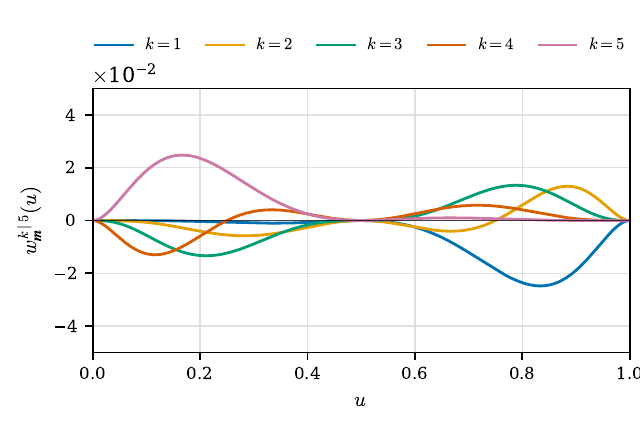}
		\caption{$\boldsymbol{m}_2=(2,2,0,0,0)$}
		\label{fig:structural-weights-22000-n5}
	\end{subfigure}\hfill
	\begin{subfigure}[t]{0.33\textwidth}
		\centering
		\includegraphics[width=\linewidth]{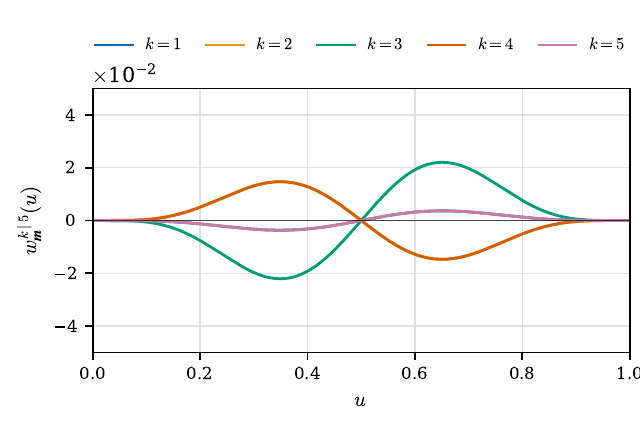}
		\caption{$\boldsymbol{m}_3=(2,1,1,1,1)$}
		\label{fig:structural-weights-21111-n5}
	\end{subfigure}
	
	\caption{Structural weights $w_{\boldsymbol{m}}^{k\mid 5}(u)$, $k=1,\ldots,5$, for three representative multi-indices.}
	\label{fig:selected-structural-weights-n5-three}
\end{figure}

\begin{example}\label{ex:selected-structural-weights-n5}
	We consider
	$	\boldsymbol{m}_1=(1,1,0,0,0)$,
	$\boldsymbol{m}_2=(2,2,0,0,0)$, and
	$\boldsymbol{m}_3=(2,1,1,1,1)$.
	The first multi-index
	is a pairwise descriptor whose
	structural weights were computed explicitly in Example~\ref{ex:structural-weights-11000}. The second one is another pairwise descriptor, but now the two active coordinates have Legendre degree \(2\). The third one is a fully active descriptor, since it involves all five components.	
	Using
	\[
	A_1^{-}(u)=-\sqrt{3}\,u(1-u),
	\qquad\textup{and}\qquad
	A_2^{-}(u)=\sqrt{5}\,u(1-u)(1-2u),
	\]
	the products of primitives appearing in \eqref{eq:sec4-weight-reduction} are
	\begin{eqnarray*}
		\prod_{i:m_{1,i}>0}A_{m_{1,i}}^{-}(u)&	=&	3u^2(1-u)^2, \\
		\prod_{i:m_{2,i}>0}A_{m_{2,i}}^{-}(u) &	=&	5u^2(1-u)^2(1-2u)^2,\\
		\prod_{i:m_{3,i}>0}A_{m_{3,i}}^{-}(u) &	=& 9\sqrt{5}\,u^5(1-u)^5(1-2u).
	\end{eqnarray*}
	Thus, the first two descriptors have \(s(\boldsymbol{m}_i)=2\), for $i=1,2$, and share the same system factor \(B_2^{k\mid 5}(u)\), whereas the third one has \(s(\boldsymbol{m}_3)=5\) and uses the fully active factor \(B_5^{k\mid 5}(u)\).
	
	Figure~\ref{fig:selected-structural-weights-n5-three} displays the signed structural weights \(w_{\boldsymbol{m}}^{k\mid 5}(u)\) for \(k=1,\ldots,5\). These curves are independent of the copula and describe how the \(k\)-out-of-\(5\) architecture transmits each dependence descriptor into the reliability expansion.
	
	The factorisation in Theorem~\ref{thm:structural-weight-factorisation} separates two sources of sign variation: the descriptor-side product of primitives and the system-side factor \(B_s^{k\mid 5}(u)\).
	For example, Panel~\ref{fig:structural-weights-11000-n5}, corresponding to the pairwise degree-one descriptor \(\boldsymbol{m}_1=(1,1,0,0,0)\), shows that the extreme architectures exhibit opposite signs and are concentrated in opposite threshold regions. Specifically, the weight of the $1$-out-of-$5$ configuration is concentrated at high thresholds, whereas the weight of the $5$-out-of-$5$ configuration is concentrated at low thresholds. This same complementary behavior is observed in the intermediate architectures,
	$k=2$ and $k=4$.
	
	Panel~\ref{fig:structural-weights-22000-n5} corresponds to the descriptor
	\(\boldsymbol{m}_2=(2,2,0,0,0)\). In this case, the product of primitives contains the factor
	\((1-2u)^2\), so all associated structural weights vanish at \(u=1/2\).
	This central cancellation follows directly from the factorisation in
	\eqref{eq:sec4-weight-reduction} and is therefore a property of the selected descriptor,
	not of the copula. It means that, at this threshold, the corresponding dependence component does not contribute to \(\mathcal{D}_{k\mid 5}(u)\).	
	However, this zero does not by itself determine whether the weight changes sign at the central threshold, because the sign also depends on the structural factor \(B_2^{k\mid 5}(u)\). Indeed, from \eqref{eq:Bkns-definition},
	\[
	B_2^{3\mid 5}(u)=3u(1-u)(2u-1),
	\]
	whereas \(B_2^{k\mid 5}(1/2)\neq0\) for \(k=1,2,4,5\). Thus, for the middle architecture, the system-side factor itself changes sign at \(u=1/2\).
	
	This is also reflected in
	the duality relation in
	Proposition~\ref{prop:dual_weights}, since \(\deg(\boldsymbol{m}_2)=4\),
	we have
	\begin{equation*}\label{eq:example-duality-22000-central}
		w_{\boldsymbol{m}_2}^{3\mid 5}(1-u)
		=
		-
		w_{\boldsymbol{m}_2}^{3\mid 5}(u),
	\end{equation*}
	so the
	structural weight of the middle system \(k=3\) is antisymmetric around \(u=1/2\), and therefore changes sign at the central threshold. 
	For the remaining
	architectures, duality pairs different systems,
	\begin{equation*}\label{eq:example-duality-22000-pairs}
		w_{\boldsymbol{m}_2}^{1\mid 5}(1-u)
		=
		-
		w_{\boldsymbol{m}_2}^{5\mid 5}(u),
		\qquad
		w_{\boldsymbol{m}_2}^{2\mid 5}(1-u)
		=
		-
		w_{\boldsymbol{m}_2}^{4\mid 5}(u),
	\end{equation*}
	and therefore this duality does not force each individual weight to change sign at
	the central threshold.
	
	Finally, Panel~\ref{fig:structural-weights-21111-n5} shows the fully active descriptor \(\boldsymbol{m}_3=(2,1,1,1,1)\). Since \(s(\boldsymbol{m}_3)=5\), all components are involved in the same tensor-product term. Here the product of primitives contains only one factor \((1-2u)\). Hence the primitive part itself changes sign at \(u=1/2\), which explains the sign reversal observed in the panel. Moreover, \(n=5\) is odd, and Proposition~\ref{prop:fully-active-dual-weights} implies that dual architectures have the same signed structural weights. Hence only three distinct curves appear, namely, those for \(k=1\) and \(k=5\) coincide, those for \(k=2\) and \(k=4\) coincide, and the \(3\)-out-of-\(5\) curve is self-dual. Thus, Panels~\ref{fig:structural-weights-22000-n5} and~\ref{fig:structural-weights-21111-n5} separate two different mechanisms: for \(\boldsymbol{m}_2\), the sign change of the middle architecture comes from \(B_2^{3\mid 5}(u)\), whereas for \(\boldsymbol{m}_3\), it comes from the primitive factor.
\end{example}

\section{Controlled reliability study under multivariate FGM dependence}
\label{sec:fgm-applications}

We next illustrate how the orthogonal expansion separates pairwise and higher-order dependence contributions to $k$-out-of-$n$ reliability. This distinction is especially relevant in redundant systems, whose reliability depends on joint component behaviour and may therefore be affected by dependence features that are not captured by pairwise summaries alone.

The multivariate FGM family provides a convenient controlled framework for this purpose, since its parameters can be associated with interactions involving different subsets of components. Combined with the representation
\[
R_{k\mid n}(u)
=
R^{\mathrm{ind}}_{k\mid n}(u)
+
\sum_{\substack{\boldsymbol{m}\in\mathbb{N}_0^n\\
		s(\boldsymbol{m})\geq 2}}
\rho_{\boldsymbol{m}}
w_{\boldsymbol{m}}^{k\mid n}(u).
\]
this allows us to isolate changes in higher-order dependence while preserving pairwise dependence and to quantify their effect on system reliability through the corresponding structural weights.

\subsection{FGM copulas in the orthogonal expansion}
\label{subsec:fgm-legendre-reliability}

Consider the standard multivariate FGM copula
\begin{equation}
	\label{eq:num-standard-multivariate-FGM}
	\mathcal{C}_{\mathrm{FGM}}(\mathbf{v})
	=
	\left(\prod_{i=1}^{n}v_i\right)
	\left[
	1+
	\sum_{\substack{A\subseteq\{1,\ldots,n\}\\ |A|\geq2}}
	\theta_A\prod_{i\in A}(1-v_i)
	\right],
\end{equation}
where the parameters satisfy
\begin{equation}
	\label{eq:num-standard-multivariate-FGM-admissibility}
	1+
	\sum_{\substack{A\subseteq\{1,\ldots,n\}\\ |A|\geq2}}
	\theta_A\prod_{i\in A}\varepsilon_i
	\geq0,
	\qquad
	(\varepsilon_1,\ldots,\varepsilon_n)\in\{-1,1\}^{n}.
\end{equation}
Each parameter $\theta_A$ is associated with the dependence among the
components indexed by the subset $A\subseteq\{1,\ldots,n\}$. The subset structure therefore provides a
natural hierarchical decomposition by dimension, namely parameters with $|A|=2$
determine the bivariate dependence, whereas parameters with $|A|\geq 3$
describe higher-dimensional dependence contributions,
see \cite{johnsonKotz1975fgm,otaKimura2021multivariateFGM}.

Differentiating \eqref{eq:num-standard-multivariate-FGM} with respect to all coordinates
gives the copula density
\begin{equation}
	\label{eq:num-standard-multivariate-FGM-density}
	c_{\mathrm{FGM}}(\mathbf{v})
	=
	1+
	\sum_{\substack{A\subseteq\{1,\ldots,n\}\\ |A|\geq2}}
	\theta_A\prod_{i\in A}(1-2v_i).
\end{equation}
By \eqref{eq:bj-def}, $\phi_1(v)=\sqrt{3}(2v-1)$, and hence every term in
\eqref{eq:num-standard-multivariate-FGM-density} is exactly a tensor-product element
of the orthonormal basis introduced in \eqref{eq:balpha-def}. More precisely, for every
$\boldsymbol{m}\in\{0,1\}^{n}$ with $s(\boldsymbol{m})\geq2$,
\[
\prod_{i\in\mathcal{A}(\boldsymbol{m})}(1-2v_i)
=
(-1)^{s(\boldsymbol{m})}
3^{-s(\boldsymbol{m})/2}
\Phi_{\boldsymbol{m}}(\mathbf{v}).
\]
Therefore, \eqref{eq:num-standard-multivariate-FGM-density} has the exact finite
orthogonal expansion
\begin{equation}
	\label{eq:num-standard-FGM-Legendre-expansion}
	c_{\mathrm{FGM}}(\mathbf{v})
	=
	1+
	\sum_{\substack{\boldsymbol{m}\in\{0,1\}^{n}\\
			s(\boldsymbol{m})\geq2}}
	(-1)^{s(\boldsymbol{m})}
	3^{-s(\boldsymbol{m})/2}
	\theta_{\mathcal{A}(\boldsymbol{m})}
	\Phi_{\boldsymbol{m}}(\mathbf{v}).
\end{equation}
Comparison with \eqref{eq:sec3-restricted-expansion} yields
\begin{equation}
	\label{eq:num-standard-FGM-Legendre-coefficients}
	\rho_{\boldsymbol{m}}
	=
	(-1)^{s(\boldsymbol{m})}
	3^{-s(\boldsymbol{m})/2}
	\theta_{\mathcal{A}(\boldsymbol{m})},
	\qquad
	\boldsymbol{m}\in\{0,1\}^{n},
	\quad
	s(\boldsymbol{m})\geq2,
\end{equation}
while all remaining nonconstant copula coefficients vanish. Thus, the subset structure
of the multivariate FGM family coincides exactly with the interaction-order structure
defined by $s(\boldsymbol{m})$. More precisely, each parameter $\theta_A$ corresponds
to a unique orthogonal coefficient indexed by a multi-index $\boldsymbol{m}$ satisfying
\[
\mathcal{A}(\boldsymbol{m})=A,
\qquad
s(\boldsymbol{m})=|A|.
\]
Finally, Theorem~\ref{prop:sec4-linear-form} gives
\begin{equation*}
	\label{eq:num-standard-FGM-reliability}
	R_{k\mid n}(u)
	=
	R_{k\mid n}^{\mathrm{ind}}(u)
	+
	\sum_{\substack{\boldsymbol{m}\in\{0,1\}^{n}\\
			s(\boldsymbol{m})\geq2}}
	(-1)^{s(\boldsymbol{m})}
	3^{-s(\boldsymbol{m})/2}
	\theta_{\mathcal{A}(\boldsymbol{m})}
	w_{\boldsymbol{m}}^{k\mid n}(u).
\end{equation*}
Hence, the FGM parameters determine the dependence contributions, while the structural
weights determine how interactions of each order are transmitted to the reliability of the
$k$-out-of-$n$ system. The finite and exact correspondence in
\eqref{eq:num-standard-FGM-Legendre-coefficients} makes the multivariate FGM family
particularly suitable for controlled comparisons of pairwise and higher-order dependence
effects.

\subsection{Same pairwise dependence, different system reliability}
\label{subsec:num-same-bivariate-margins}

We consider the five $k$-out-of-$5$ system architectures, for
$k=1,\ldots,5$, under two symmetric FGM dependence specifications,
referred to as Specification~I and Specification~II. For any multi-index
$\boldsymbol{m}$ with $m_i=m_j=1$ and all remaining entries equal to
zero, the coefficient $\rho_{\boldsymbol{m}}$ coincides with the
Spearman correlation of the corresponding bivariate margin
$(U_i,U_j)$; see \cite{ngounou2025orthproj}. Since every active coordinate in the standard FGM expansion
\eqref{eq:num-standard-FGM-Legendre-expansion} is represented by the
degree-one univariate orthonormal polynomial $\phi_1$, the coefficients
with $s(\boldsymbol{m})=2$ are precisely the pairwise Spearman
correlations.

We define Specification~I by assigning the same coefficient to all
terms with $s(\boldsymbol{m})=2$ and setting all higher-order
coefficients equal to zero. By symmetry,
$\boldsymbol{m}_{(I)}=(1,1,0,0,0)$ can be taken as a representative
pairwise multi-index. Hence, up to permutation of the component labels,
all component pairs have the same Spearman correlation,
\[
\rho_S(U_i,U_j)
=
\rho_{\boldsymbol{m}_{(I)}},
\qquad i\neq j.
\]
From \eqref{eq:num-standard-FGM-Legendre-coefficients}, the common
pairwise coefficient corresponds to
\begin{equation}
	\label{eq:num-same-pairwise-theta-I}
	\theta_A
	=
	3\rho_{\boldsymbol{m}_{(I)}},
	\qquad |A|=2.
\end{equation}
Substitution into \eqref{eq:num-standard-FGM-Legendre-expansion} gives
\begin{equation*}
	\label{eq:num-same-pairwise-density-I}
	c_{\mathrm{FGM}}^{(I)}(\mathbf{v})
	=
	1+
	\rho_{\boldsymbol{m}_{(I)}}
	\sum_{\substack{\boldsymbol{m}\in\{0,1\}^{5}\\
			s(\boldsymbol{m})=2}}
	\Phi_{\boldsymbol{m}}(\mathbf{v}).
\end{equation*}

We now consider Specification~II, which preserves the pairwise
dependence structure of Specification~I and adds a common
four-component dependence contribution. By symmetry,
$\boldsymbol{m}_{(II)}=(1,1,1,1,0)$ can be taken as a representative
multi-index for the terms with $s(\boldsymbol{m})=4$. From
\eqref{eq:num-standard-FGM-Legendre-coefficients}, the common
four-component coefficient corresponds to
\begin{equation}
	\label{eq:num-same-pairwise-theta-II}
	\theta_A
	=
	9\rho_{\boldsymbol{m}_{(II)}},
	\qquad |A|=4.
\end{equation}
Substitution into \eqref{eq:num-standard-FGM-Legendre-expansion} gives
\begin{equation}
	\label{eq:num-same-pairwise-density-II}
	c_{\mathrm{FGM}}^{(II)}(\mathbf{v})
	=
	c_{\mathrm{FGM}}^{(I)}(\mathbf{v})
	+
	\rho_{\boldsymbol{m}_{(II)}}
	\sum_{\substack{\boldsymbol{m}\in\{0,1\}^{5}\\
			s(\boldsymbol{m})=4}}
	\Phi_{\boldsymbol{m}}(\mathbf{v}).
\end{equation}

For the standard multivariate FGM family, bivariate margins depend only
on the corresponding pairwise parameters; see, e.g.,
\cite{otaKimura2021multivariateFGM}. Hence, since Specifications~I and~II
share all parameters with $|A|=2$, they have identical bivariate
copulas. More explicitly, for every pair $i\neq j$,
\begin{equation*}
	\label{eq:num-same-pairwise-common-bivariate}
	\mathcal{C}_{ij}^{(I)}(v_i,v_j)
	=
	\mathcal{C}_{ij}^{(II)}(v_i,v_j)
	=
	v_i v_j
	\left[
	1+
	3\rho_{\boldsymbol{m}_{(I)}}
	(1-v_i)(1-v_j)
	\right].
\end{equation*}
Thus, the complete bivariate copula is held fixed across the two specifications.

The parameters of the copula density
\eqref{eq:num-same-pairwise-density-II} must satisfy the FGM
admissibility condition in
\eqref{eq:num-standard-multivariate-FGM-admissibility}. Using
\eqref{eq:num-same-pairwise-theta-I} and
\eqref{eq:num-same-pairwise-theta-II}, this condition becomes
\[
1
+
3\rho_{\boldsymbol{m}_{(I)}}
\sum_{\substack{A\subseteq\{1,\ldots,5\}\\ |A|=2}}
\prod_{i\in A}\varepsilon_i
+
9\rho_{\boldsymbol{m}_{(II)}}
\sum_{\substack{A\subseteq\{1,\ldots,5\}\\ |A|=4}}
\prod_{i\in A}\varepsilon_i
\geq 0,
\]
for every $(\varepsilon_1,\ldots,\varepsilon_5)\in\{-1,1\}^5$.
Because both specifications are symmetric with respect to the component
labels, the $32$ sign configurations reduce to three distinct
admissibility conditions:
\begin{equation*}
	\label{eq:num-same-pairwise-admissibility}
	\begin{aligned}
		1
		+30\rho_{\boldsymbol{m}_{(I)}}
		+45\rho_{\boldsymbol{m}_{(II)}} &\geq 0,\\
		1
		+6\rho_{\boldsymbol{m}_{(I)}}
		-27\rho_{\boldsymbol{m}_{(II)}} &\geq 0,\\
		1
		-6\rho_{\boldsymbol{m}_{(I)}}
		+9\rho_{\boldsymbol{m}_{(II)}} &\geq 0.
	\end{aligned}
\end{equation*}

By Remark~\ref{rem:dependence-impact-differences}, the reliability
difference between the two specifications is also their difference in
dependence impact. Since Specification~II differs from
Specification~I only through the five terms with
$s(\boldsymbol{m})=4$, \eqref{eq:dependence-impact-coefficient-difference}
and the structural-weight invariance in
\eqref{eq:structural-weight-permutation-invariance} give
\begin{equation}
	\label{eq:num-same-pairwise-reliability-difference}
	R_{k\mid5}^{(II)}(u)
	-
	R_{k\mid5}^{(I)}(u)
	=
	D_{k\mid5}^{(II)}(u)
	-
	D_{k\mid5}^{(I)}(u)
	=
	5\rho_{\boldsymbol{m}_{(II)}}
	w_{\boldsymbol{m}_{(II)}}^{k\mid5}(u).
\end{equation}
Thus, after fixing the complete bivariate dependence structure, the
change in reliability is exactly the change in dependence impact
generated by the four-component contribution.

For $\boldsymbol{m}_{(II)}=(1,1,1,1,0)$,
Theorem~\ref{thm:structural-weight-factorisation}, together with
\eqref{eq:sec4-weight-reduction} and
\eqref{eq:Bkns-definition}, yields
\begin{equation*}
	\label{eq:num-same-pairwise-explicit-differences}
	\begin{aligned}
		R_{1\mid5}^{(II)}(u)-R_{1\mid5}^{(I)}(u)
		&=
		-45\rho_{\boldsymbol{m}_{(II)}}
		u^{5}(1-u)^{4},
		\\
		R_{2\mid5}^{(II)}(u)-R_{2\mid5}^{(I)}(u)
		&=
		45\rho_{\boldsymbol{m}_{(II)}}
		u^{4}(1-u)^{4}(4u-1),
		\\
		R_{3\mid5}^{(II)}(u)-R_{3\mid5}^{(I)}(u)
		&=
		135\rho_{\boldsymbol{m}_{(II)}}
		u^{4}(1-u)^{4}(1-2u),
		\\
		R_{4\mid5}^{(II)}(u)-R_{4\mid5}^{(I)}(u)
		&=
		45\rho_{\boldsymbol{m}_{(II)}}
		u^{4}(1-u)^{4}(4u-3),
		\\
		R_{5\mid5}^{(II)}(u)-R_{5\mid5}^{(I)}(u)
		&=
		45\rho_{\boldsymbol{m}_{(II)}}
		u^{4}(1-u)^{5}.
	\end{aligned}
\end{equation*}

For $\rho_{\boldsymbol{m}_{(II)}}>0$, Specification~II yields lower
reliability than Specification~I for the $1$-out-of-$5$ architecture
throughout the full operating range on the standardised lifetime scale,
whereas the opposite relation holds for the $5$-out-of-$5$ architecture.
For the intermediate architectures, the sign of
$R_{k\mid5}^{(II)}(u)-R_{k\mid5}^{(I)}(u)$ changes at $u=1/4$,
$u=1/2$, and $u=3/4$ for $k=2$, $k=3$, and $k=4$, respectively.
Moreover, all five differences contain the factor $u^4(1-u)^4$,
which is consistent with Proposition~\ref{prop:wm-extremes} and, in
particular, with the boundary bounds in
\eqref{eq:wm-extremes-lower-upper} for
$s(\boldsymbol{m}_{(II)})=4$.

Equation~\eqref{eq:num-same-pairwise-reliability-difference} also provides
a direct measure of the largest reliability discrepancy between the two
dependence specifications over the entire operating range. Namely,
\begin{equation}
	\label{eq:num-same-pairwise-maximum-reliability-difference}
	\max_{u\in(0,1)}
	\left|
	R_{k\mid5}^{(II)}(u)
	-
	R_{k\mid5}^{(I)}(u)
	\right|
	=
	5|\rho_{\boldsymbol{m}_{(II)}}|
	\max_{u\in(0,1)}
	|	w_{\boldsymbol{m}_{(II)}}^{k\mid5}(u)	|.
\end{equation}
Since the complete bivariate dependence structure is fixed, the quantity
in \eqref{eq:num-same-pairwise-maximum-reliability-difference} measures
the maximum change in system reliability attributable to the
four-component dependence contribution. It also provides a direct
measure of the sensitivity of each $k$-out-of-$5$ architecture to this
higher-order dependence component.

\begin{figure}[ht]
	\centering
	\includegraphics[width=0.5\textwidth]{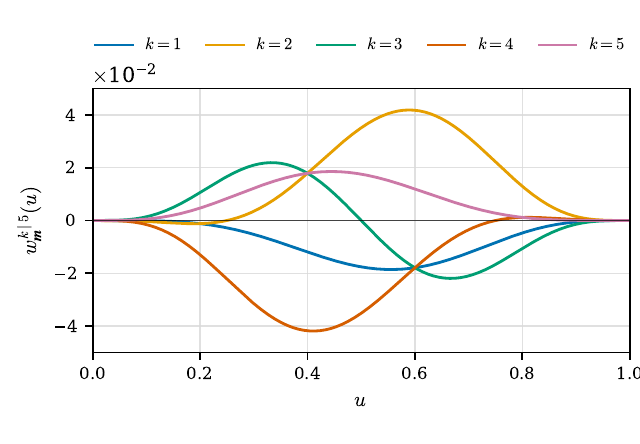}
	\caption{Structural weights $w_{\boldsymbol{m}_{(II)}}^{k\mid5}(u)$ associated with the four-component descriptor $\boldsymbol{m}_{(II)}=(1,1,1,1,0)$, plotted for the five $k$-out-of-$5$ architectures, $k=1,\ldots,5$.}
	\label{fig:structural-weights-mII}
\end{figure}

In view of
\eqref{eq:num-same-pairwise-maximum-reliability-difference},
Table~\ref{tab:max-structural-weights-mII} reports the maximum absolute
structural weight and the threshold values at which it is attained
for each $k$-out-of-$5$ architecture, while
Figure~\ref{fig:structural-weights-mII} displays the corresponding
profiles over the complete operating range.

\begin{table}[ht]
	\centering
	\caption{Maximum absolute structural weights for
		$\boldsymbol{m}_{(II)}=(1,1,1,1,0)$ across the five
		$k$-out-of-$5$ architectures.}
	\label{tab:max-structural-weights-mII}
	\begin{tabular}{c c c}
		\hline
		Architecture $k$
		&
		Maximising threshold(s) $u$
		&
		$\displaystyle \max_{u\in(0,1)}
		|w_{\boldsymbol{m}_{(II)}}^{k\mid5}(u)|$
		\\
		\hline
		$1$ & $5/9$ & $0.01858$ \\
		$2$ & $(7+\sqrt{13})/18$ & $0.04191$ \\
		$3$ & $1/3,\ 2/3$ & $0.02195$ \\
		$4$ & $(11-\sqrt{13})/18$ & $0.04191$ \\
		$5$ & $4/9$ & $0.01858$ \\
		\hline
	\end{tabular}
\end{table}

The reflected structure of the maximising thresholds reported in
Table~\ref{tab:max-structural-weights-mII} follows directly from the
structural duality established in Proposition~\ref{prop:dual_weights}.
Since $\deg(\boldsymbol{m}_{(II)})=4$, \eqref{eq:dual_weight_relation}
specialises, for $n=5$, to
\begin{equation*}
	\label{eq:num-mII-duality}
	w_{\boldsymbol{m}_{(II)}}^{k\mid5}(1-u)
	=
	-w_{\boldsymbol{m}_{(II)}}^{6-k\mid5}(u).
\end{equation*}
In particular, $k=3$ is self-dual, and therefore
\begin{equation*}
	\label{eq:num-mII-self-dual-antisymmetry}
	w_{\boldsymbol{m}_{(II)}}^{3\mid5}(1-u)
	=
	-w_{\boldsymbol{m}_{(II)}}^{3\mid5}(u).
\end{equation*}
Hence,
$|w_{\boldsymbol{m}_{(II)}}^{3\mid5}(u)|$ is symmetric
about $u=1/2$, consistently with
Corollary~\ref{cor:dual_weight_magnitude_symmetry}, and the two equal
absolute maxima are attained at the reflected thresholds
$u=1/3$ and $u=2/3$. The same duality explains the remaining pairs
reported in Table~\ref{tab:max-structural-weights-mII}. In particular,
the $1$-out-of-$5$ and $5$-out-of-$5$ architectures have the same
maximum absolute weight at the reflected thresholds $u=5/9$ and
$u=4/9$, respectively.
Likewise, the $2$-out-of-$5$ and $4$-out-of-$5$ architectures attain
their common maximum absolute weight at the reflected thresholds
$u=(7+\sqrt{13})/18$ and $u=(11-\sqrt{13})/18$, respectively.

The comparison reveals a clear architecture-dependent response to the
four-component dependence contribution. The $2$-out-of-$5$ and
$4$-out-of-$5$ architectures attain the largest absolute structural
weight,
approximately $0.04191$, followed by the $3$-out-of-$5$ architecture,
with approximately $0.02195$, whereas the $1$-out-of-$5$ and
$5$-out-of-$5$ architectures attain the smallest maximum,
approximately $0.01858$. Hence, for a fixed value of
$|\rho_{\boldsymbol{m}_{(II)}}|$, the $2$-out-of-$5$ and
$4$-out-of-$5$ systems are the most sensitive among the five
architectures to the additional four-component dependence considered
here.

The results therefore show that the four-component dependence
contribution is transmitted very differently across the five
$k$-out-of-$5$ architectures. Its effect is strongest for the
$2$-out-of-$5$ and $4$-out-of-$5$ systems, changes sign at
architecture specific threshold levels, and is attenuated near the
boundaries of the operating range. These patterns are consistent
with the structural weight properties established in
Section~\ref{sec:linear-form-weights}, including boundary attenuation
and duality.


\section{Conclusions}
\label{sec:conclusions}

This research has developed an orthogonal representation of the reliability
of a $k$-out-of-$n$ system with dependent components and identical
continuous marginal lifetime distributions. Under an absolutely
continuous copula with square-integrable density, the standardised
reliability is expressed as the independence benchmark plus a linear
combination of copula coefficients specified in  \eqref{eq:sec4-reliability-linear-form}.
This representation separates the dependence information contained in
the coefficients $\rho_{\boldsymbol{m}}$ from the response of the
system architecture, which is entirely described by the structural
weights $w_{\boldsymbol{m}}^{k\mid n}(u)$. Consequently, the influence
of a given dependence component on reliability is determined jointly by
its magnitude in the copula expansion and by the way in which the
$k$-out-of-$n$ structure transmits that component at the operating
threshold.

The dependence impact function introduced in
Definition~\ref{def:dependence-impact} provides a system-level measure
of the departure from independence. Besides its pointwise
interpretation, Proposition~\ref{prop:mean-life-impact} links its
integral on the original lifetime scale to the corresponding change in
mean system lifetime. Proposition~\ref{prop:sum-D-zero} further shows
that, at any fixed threshold, dependence redistributes its reliability
impact across the family of $k$-out-of-$n$ architectures. Exchangeability
permits the orthogonal expansion to be regrouped over permutation
classes, whereas radial symmetry relates the dependence impacts of dual
architectures through threshold reflection.

A central result is the factorisation established in
Theorem~\ref{thm:structural-weight-factorisation}. Each structural
weight separates into a product of one-dimensional Legendre primitives
associated with the dependence descriptor and the factor
$B_{s(\boldsymbol{m})}^{k\mid n}(u)$ associated with the system
architecture. This factorisation avoids direct multivariate integration
and makes several structural properties explicit. In particular,
Proposition~\ref{prop:wm-extremes} shows that the contribution of a
descriptor involving $s(\boldsymbol{m})$ components is attenuated at
least at orders $u^{s(\boldsymbol{m})}$ and
$(1-u)^{s(\boldsymbol{m})}$ in the extreme reliability regimes, corresponding respectively to high-survival and
low-survival operating conditions. The duality
relations in Proposition~\ref{prop:dual_weights} and
Corollary~\ref{cor:dual_weight_magnitude_symmetry} determine the
response of dual architectures at reflected thresholds, while
Proposition~\ref{prop:fully-active-dual-weights} gives a stronger
pointwise relation when the descriptor involves all components.

The multivariate FGM study illustrates how the orthogonal representation
provides a quantitative system-level characterisation of higher-order
dependence beyond the information contained in the bivariate margins.
For the FGM construction considered here, the four-component dependence
contribution is isolated explicitly and its effect on reliability is
transmitted exactly through
\eqref{eq:num-same-pairwise-reliability-difference}. Consequently, the
structural weight
$w_{\boldsymbol{m}_{(II)}}^{k\mid5}(u)$ determines, for each system
architecture and operating threshold, both the magnitude and the
direction of the resulting reliability change. The attenuation of this
effect in the extreme reliability regimes and the reflected behaviour
of dual architectures are explained by the structural properties
established in Proposition~\ref{prop:wm-extremes},
Proposition~\ref{prop:dual_weights}, and
Corollary~\ref{cor:dual_weight_magnitude_symmetry}. The application thus
shows how the theoretical properties of the structural weights can be
used to characterise quantitatively the system-level transmission of a
specified higher-order dependence component.
The \textit{same pairwise dependence} setting considered here is therefore one
informative instance of a more general framework, since the proposed
reliability methodology applies to arbitrary dependence structures whose
copula densities admit the orthogonal representation recalled in
Subsection~\ref{subsec:legendre-notation}.

The present analysis is restricted to identical continuous margins and
absolutely continuous copulas with square-integrable densities. Natural
extensions include systems with non-identical marginal lifetime
distributions and more general coherent structures. From a statistical
perspective, an important further step is to combine the structural
representation with data-driven estimation of the copula coefficients,
together with uncertainty quantification for the induced reliability
correction. Such developments would allow the orthogonal decomposition
to be used not only as an analytical representation of dependence
effects, but also as a practical tool for identifying which estimated
dependence components are relevant for system-level reliability.


\section*{Funding}

This research was supported by Grant CEX2023-001347-S,
funded by \seqsplit{MICIU/AEI/10.13039/501100011033}.

\section*{Declaration of competing interests}

The author declares no competing interests.


\end{document}